\documentclass{article}
\usepackage{amssymb}
\usepackage{amsmath,color}
\usepackage[hidelinks]{hyperref}
\usepackage{easyReview}
\usepackage{cancel}

\usepackage[square,numbers]{natbib}         
\usepackage{ntheorem}
\newtheorem{lem}{Lemma}[section]
\newtheorem{cor}[lem]{Corollary}
\newtheorem{teo}[lem]{Theorem}
\newtheorem*{teoNon}{Theorem}
\newtheorem{os}[lem]{Remark}

\newtheorem{prop}[lem]{Proposition}

\newcommand{\qed}{\thinspace\null\nobreak\hfill\hbox{\vbox{\kern-.2pt\hrule
 height.2pt depth.2pt\kern-.2pt\kern-.2pt \hbox to2.5mm{\kern-.2pt\vrule
 width.4pt \kern-.2pt\raise2.5mm\vbox to.2pt{}\lower0pt\vtop
 to.2pt{}\hfil\kern-.2pt \vrule
 width.4pt \kern-.2pt}\kern-.2pt\kern-.2pt\hrule height.2pt depth.2pt
 \kern-.2pt}}\par\medbreak}

\newcommand{\R}{\mathbb{R}}

\newcommand{\C}{\mathbb{C}}

\newcommand{\Rp}{\textrm{\emph{Re}\,}}
\newcommand{\Ip}{\textrm{\emph{Im}\,}}
\newcommand{\eps}{\varepsilon}

\newcommand{\ds}{\displaystyle}

\date{}
\begin{document}

\title{Regularity theory for   parabolic operators in the half-space with  boundary degeneracy}
\author{G. Metafune \thanks{Dipartimento di Matematica e Fisica ``Ennio De Giorgi'', Universit\`a del Salento, C.P.193, 73100, Lecce, Italy.
-mail:  giorgio.metafune@unisalento.it}\qquad L. Negro \thanks{Dipartimento di Matematica e Fisica  ``Ennio De
Giorgi'', Universit\`a del Salento, C.P.193, 73100, Lecce, Italy. email: luigi.negro@unisalento.it} \qquad C. Spina \thanks{Dipartimento di Matematica e Fisica``Ennio De Giorgi'', Universit\`a del Salento, C.P.193, 73100, Lecce, Italy.
e-mail:  chiara.spina@unisalento.it}}

\maketitle
\begin{abstract}
\noindent 
We study elliptic and parabolic problems governed by the  singular elliptic   operators 
\begin{align*}
	\mathcal L=y^{\alpha_1}\mbox{Tr }\left(QD^2_xu\right)+2y^{\frac{\alpha_1+\alpha_2}{2}}q\cdot \nabla_xD_y+\gamma y^{\alpha_2} D_{yy}+{cy^{\alpha_2-1}D_y}
\end{align*}
under {Neumann}   boundary condition,
in the half-space $\R^{N+1}_+=\{(x,y): x \in \R^N, y>0\}$. We prove elliptic and
parabolic $L^p$-estimates and solvability for the associated problems. In the language of semigroup
theory, we prove that $\mathcal L$ generates an analytic semigroup, characterize its domain as a weighted
Sobolev space and show that it has maximal regularity.

\par

\end{abstract}

\section{Introduction}

In this paper we study solvability and regularity of elliptic and parabolic problems associated to the  degenerate   operators
\begin{align}\label{defL}
	\mathcal L=y^{\alpha_1}\mbox{Tr }\left(QD^2_xu\right)+2y^{\frac{\alpha_1+\alpha_2}{2}}q\cdot \nabla_xD_y+\gamma y^{\alpha_2}D_{yy}+{cy^{\alpha_2-1}D_y}
\end{align}
and $D_t- \mathcal L$ in the half-space $\R^{N+1}_+=\{(x,y): x \in \R^N, y>0\}$ or  in $(0, \infty) \times \R^{N+1}_+$ and under {Neumann} boundary condition at $y=0$.

Here {$c\in\R$}  and 
$
\left(
\begin{array}{c|c}
	Q  & { q}^t \\[1ex] \hline
	q& \gamma
\end{array}\right)$ is 
a constant real  elliptic matrix. The real numbers $\alpha_1, \alpha_2$ satisfy $\alpha_2<2$ and $\alpha_2-\alpha_1<2$ but are not assumed to be nonnegative.  

We write
$
B_y
$ to denote the 1-dimensional  Bessel operator $D_{yy}+\frac{c}{y}D_y$. With this notation  the   special case where
\begin{equation*} 
\mathcal L=y^{\alpha_1}\Delta_x+y^{\alpha_2}B_y
\end{equation*}
has been already studied in \cite{MNS-CompleteDegenerate}.  The main novelty here consists in the presence of the mixed derivatives $2y^{\frac{\alpha_1+\alpha_2}{2}}q\cdot \nabla_xD_y$ in the operator $\mathcal L$ which is a crucial step for treating degenerate operators in domains, through a localization procedure.

\medskip
Our main result is the following,   
see Theorems \ref{complete-Bessel}, \ref {complete-oblique} and  Appendix B for the definition of the weighted Sobolev spaces involved.
\begin{teoNon} 
	Let $\alpha_1,\alpha_2\in\R$ such that $\alpha_2<2$,  $\alpha_2-\alpha_1<2$ and
	$$\alpha_1^{-} <\frac{m+1}p<\frac{c}{\gamma}+1-\alpha_2.$$   
	Then the operator
	\begin{align*}
		\mathcal L=y^{\alpha_1}\mbox{Tr }\left(Q D^2_xu\right)+2y^{\frac{\alpha_1+\alpha_2}{2}}q\cdot \nabla_xD_y+\gamma y^{\alpha_2} D_{yy}+{c y^{\alpha_2-1}D_y}
	\end{align*} endowed with domain {$W^{2,p}_{\mathcal N}\left(\alpha_1,\alpha_2,m\right)$}   generates a bounded analytic semigroup  in $L^p_m$ which has maximal regularity.
\end{teoNon}
Let us explain the meaning of the restrictions $\alpha_2<2$, $\alpha_2-\alpha_1<2$
considering first the case where $\alpha_1=\alpha_2=\alpha$, so that the unique requirement is $\alpha<2$. 

It turns out that when $\alpha \geq 2$ the problem is easily treated in the strip $\R^N\times [0,1]$ in the case of the Lebesgue measure, see \cite{FornMetPallScn5}, and all problems are due to the strong diffusion at infinity. The case $\alpha \geq 2$ in the strip $\R^N \times [1, \infty[$ requires therefore new investigation even though the 1-dimensional case is easily treated by the change of variables of Appendix B.

When $\alpha_1 \neq \alpha_2$, the further restriction $\alpha_2-\alpha_1<2$ comes from the change of variables of Appendix B, see Section 6.

\medskip
Let us briefly describe the previous literature on these  operators. In \cite{MNS-Caffarelli, MNS-Caf-Schauder} we considered the simplest case of  $\Delta_x+B_y$ making extensive use of  the commutative structure of the operator. The non-commutative case of $y^{\alpha_1}\Delta_x+y^{\alpha_2}B_y$ have been later faced in \cite{MNS-CompleteDegenerate}.  
Another source of non-commutativity comes from the presence of mixed derivatives. In \cite{MNS-Singular-Half-Space}, we treated the operator 
$$
\mbox{Tr }\left(Q D^2_xu\right)+2q\cdot \nabla_xD_y+b \cdot \nabla_x+B_y
$$
under Neumann boundary conditions. We refer the reader also to  \cite{dong2020neumann}, \cite{dong2020RMI}, \cite{dong2020parabolic} and \cite{dong2021weighted} for related results with different methods, mainly without the power $y^\alpha$.

\medskip
This paper is devoted to a final step in this direction, by adding (different) powers of $y$ in front of the main terms of the operator.
This is, by no means, an immediate generalization of the previous  results or methods and many extra difficulties appear. Here we consider only constant matrices $Q$ and constant $q, \gamma$. The general case where $Q,q, \gamma$ are bounded and uniformly continuous is however straightforward and allows to treat  operators in smooth domains, whose degeneracy in the top order coefficients behaves like a power of the distance from the boundary. We shall treat these topics in a forthcoming paper.

%

As for the simpler operator $y^{\alpha_1}\Delta_x+y^{\alpha_2}B_y$, the case $\alpha_1=\alpha_2$ implies all other cases by the change of variables described in Appendix B. However this modifies the underlying measure and the procedure works if one is able to deal with the complete scale of $L^p_m$ spaces, where $L^p_m=L^p(\R^{N+1}_+; y^m dxdy)$.


\bigskip

The operators $\mathcal L$, $ D_t-\mathcal L$, with $\alpha_1=\alpha_2=\alpha$,   are studied through estimates like 
\begin{equation} \label{closedness}
\|y^\alpha \Delta_x u\|_{p,m}+\|y^\alpha \nabla_xD_y u\|_{p,m}+\|y^\alpha B_y u\|_{p,m} \le C\| \mathcal L u\|_{p,m}, \quad 
\end{equation}
and
\begin{equation} \label{closedness1}
\|D_t u\|_{p,m}+\|\mathcal Lu\|_{p,m} \le C\| (D_t-\mathcal L) u\|_{p,m}, \quad 
\end{equation}
where the $L^p$ norms are taken over $\R_+^{N+1}$ and on $(0, \infty) \times \R_+^{N+1}$ respectively. This kind of estimates are quite natural in this context but not easy to prove. Of course they imply $\|y^\alpha D_{x_ix_j}u\|_{p,m} \le C\| \mathcal Lu\|_{p,m}$, by the Calder\'{o}n-Zygmund inequalities in the $x$-variables. 

\bigskip
Let us explain how to obtain \eqref{closedness}. Assuming that $y^\alpha(\Delta_x u+2a\cdot\nabla_xD_y u+B_yu)=f$ and taking the Fourier transform with respect to $x$ (with covariable $\xi$) we obtain $-|\xi|^2 \hat u(\xi,y)+2ia\cdot\xi D_y+B_y \hat u(\xi,y)=y^{-\alpha}\hat f(\xi,y)$.
Denoting by ${\cal F}$ the Fourier transform with respect to $x$  we get 
\begin{align*}
& y^\alpha\Delta_x  \mathcal  L^{-1}=-{\cal F}^{-1}\left (y^\alpha|\xi|^2(|\xi|^2-2ia\xi D_y-B_y)^{-1} y^{-\alpha} \right) {\cal F}\\&
 y^\alpha\nabla_xD_y  \mathcal  L^{-1}=-{\cal F}^{-1}\left (y^\alpha\xi D_y(|\xi|^2-2ia\xi D_y-B_y)^{-1} y^{-\alpha} \right) {\cal F}
\end{align*}
and the boundedness of $y^\alpha\Delta_x  \mathcal L^{-1}$, $y^\alpha\nabla_xD_y  \mathcal L^{-1}$ are equivalent to that of the multipliers $\xi \in \R^N \to y^\alpha|\xi|^2(|\xi|^2-2i a\cdot\xi D_y -B_y)^{-1}y^{-\alpha} $ and $\xi \in \R^N \to y^\alpha \xi D_y(|\xi|^2-2ia\cdot\xi D_y -B_y)^{-1}y^{-\alpha} $ in $L^p_m$. 

The structure of these multipliers show the difficulties connected to the presence of the mixed operators. The parameter $\xi$ not only appears as a spectral parameter but also in the  operator $L_{2a\cdot \xi}:=B_y+2ia\cdot \xi D_y$. For this reason, we need a careful study of the 1- dimensional operator $L_{2a\cdot \xi}$ and  variants like $y^\alpha (L_{2a\cdot \xi}-|\xi|^2)$, with estimates which depend explicitly on $\xi$.

Both the elliptic and parabolic estimates above  share the name ``maximal regularity" even though this term is often restricted to the parabolic case. 
We refer to  \cite{KunstWeis} and the new books \cite{WeisBook1}, \cite{WeisBook2} for 
the functional analytic approach to  maximal regularity and to \cite{DenkHieberPruss} for applications of these methods to uniformly parabolic operators.

\medskip

\medskip

The paper is organized as  follows. In Section 2 we recall some results concerning a one-dimensional  Bessel operator $y^\alpha B_y$ perturbed by a potential. 
In Section 3 we define and study a  $1d$ auxiliary operator   through a quadratic form. 
In Section 4  we investigate the boundedness of some multipliers related to the degenerate operator.
In Section 5, which is the core of the paper,  we prove generation results, maximal regularity and domain characterization for the operator $\mathcal L$, under Neumann boundary conditions. 
Finally, in Section 6, we extend our results to more general operators.

In the Appendices we briefly recall the harmonic analysis background needed in the paper,  as  square function estimates, $\mathcal R$-boundedness, a vector valued multiplier theorem and the changes of variables needed to reduce our operators to the simpler case where $\alpha_1=\alpha_2$.

\bigskip
\noindent\textbf{Notation.} For $N \ge 0$, $\R^{N+1}_+=\{(x,y): x \in \R^N, y>0\}$. For $m \in \R$ we consider the measure $y^m dx dy $ in $\R^{N+1}_+$ and  we write $L^p_m(\R_+^{N+1})$ for  $L^p(\R_+^{N+1}; y^m dx dy)$ and often only $L^p_m$ when $\R^{N+1}_+$ is understood. 
$\C^+=\{ \lambda \in \C: \Rp \lambda >0 \}$ and, for $|\theta| \leq \pi$, we denote by  $\Sigma_{\theta}$  the open sector $\{\lambda \in \C: \lambda \neq 0, \ |Arg (\lambda)| <\theta\}$.  
We denote by  $\alpha^+$ and $\alpha^-$ the positive and negative  part of a real  number, that is $\alpha^+=\max\{\alpha,0\}$,  $\alpha^-=-\min\{\alpha,0\}$.

We use $B$ for the one-dimensional Bessel operator $D_{yy} +\frac cy D_y$ and $L_b$ for $B+ibD_y$. Here $c, b \in \R$ and both operators are defined on the half-line $(0, \infty)$.

\bigskip

 \section{The  1-dimensional operator $y^\alpha B-\mu y^{\alpha}$, $\mu\geq 0$}
 In this section we summarize  the main   results  proved in \cite{MNS-PerturbedBessel}  for the one dimensional operator  $y^{\alpha} B-\mu y^{\alpha}=y^\alpha\left(D_{yy}+\frac{c}{y}D_y \right)-\mu y^\alpha$, $\mu \geq 0$,  in $L^p_m$. To characterize the domain for $\mu>0$,  we note that the domain of the potential $V(y)=y^\alpha$ in $L^p_m$ is $$\left\{u\in L^p_m:\ y^{\alpha}u\in L^p_m\right\}=L^p_m \cap L^p_{m+\alpha p}.$$ We recall that the Sobolev spaces  $W^{2,p}_{\mathcal N}(\alpha,m)$  are defined in Appendix B and the pedix $\mathcal N$ indicates a Neumann boundary condition at $y=0$.
 
 \begin{teo}  \label{1d} Let  $\alpha<2$, $c\in\R$ and $1<p<\infty$.    
 	\begin{itemize}
 		\item[(i) ]If  $0<\frac{m+1}p<c+1-\alpha$, then  the operator 
 		$y^\alpha B$ endowed with domain $W^{2,p}_{\mathcal N}(\alpha,m)$	
 		generates a bounded analytic semigroup of angle $\pi/2$ on $L^p_m$ which is  positive for $z>0$. 
 		\item[(ii)] If $\mu>0$ and $\alpha^-<\frac{m+1}p<c+1-\alpha$ then  the operator 	 
 		$y^{\alpha} B-\mu y^{\alpha}$  endowed with domain $W^{2,p}_{\mathcal N}(\alpha,m)\cap  L^p_{m+\alpha p}$ generates a bounded analytic semigroup  in $L^p_m$  which is  positive for $z>0$.
 	\end{itemize}
 	In both cases the set
	\begin{equation} \label{defD}
	\mathcal{D}=\left \{u \in C_c^\infty ([0, \infty)), \ D_y u(y)=0\  {\rm for} \ y \leq \delta\ {\rm  and \ some\ } \delta>0\right \}
\end{equation}	is a core and the semigroups have maximal regularity.
 \end{teo}
 {\sc Proof.} See \cite[Theorem 4.2]{MNS-PerturbedBessel} for (i), \cite[Theorem  8.5]{MNS-PerturbedBessel} for (ii). The density of $\mathcal D$ is proved in proved in \cite[Propositions 4.3,  Theorem 8.5]{MNS-PerturbedBessel} and maximal regularity in \cite[Theorem 7.2]{MNS-PerturbedBessel}.
 \qed

The following kernel estimates for the Bessel operator $B$ (corresponding to $\alpha=\mu=0$) will be used extensively.

\begin{prop}\label{Estimates Bessel kernels}
	Let $c>-1$. The semigroup $(e^{zB})_{z \in \C_+}$ consists of integral operators.  Its heat kernel $p_B$, written  with respect the measure $\rho^cd\rho$, satisfies for every $\eps>0$, $z\in\Sigma_{\frac{\pi}{2}-\eps}$ and some  $C_\eps, \kappa_\eps>0$, 
	\begin{align*} 
		| p_{B}(z,y,\rho)|
		&\leq 
		C_\eps |z|^{-\frac{1}{2}}\rho^{-c}  \left (\frac{\rho}{|z|^{\frac{1}{2}}}\wedge 1 \right)^{c}
		\exp\left(-\frac{|y-\rho|^2}{\kappa_\eps |z|}\right),\\[1ex] 
		| D_y p_{B}(z,y,\rho)|
		&\leq 
		C_\eps |z|^{-1}\rho^{-c}\left (\frac{y}{|z|^{\frac{1}{2}}}\wedge 1 \right)\left (\frac{\rho}{|z|^{\frac{1}{2}}}\wedge 1 \right)^{c}
		\exp\left(-\frac{|y-\rho|^2}{\kappa_\eps |z|}\right).
	\end{align*}
\end{prop}
{\sc{Proof.}} See \cite[Propositions 2.8 and 2.9]{MNS-Caffarelli}. \qed

\section{The 1-dimensional operator $y^\alpha (B+2ia\cdot\xi D_y-|\xi|^2) $}\label{section L_b}
In this section we prove  generation properties in $L^p_m$ and  heat kernel bounds for the operator   $$y^\alpha (B+2ia\cdot\xi D_y-|\xi|^2)=y^\alpha L_{2a\cdot\xi}-y^\alpha|\xi|^2, \qquad \ |a|<1, \alpha<2,$$ 
where, according to our notation, $L_b=B+ibD_y$.  

Having in mind the study of the operator $y^\alpha\Delta_x u+y^\alpha2a\cdot \nabla_xD_yu+ y^\alpha B_yu$ (see  Section \ref{Sect DOm max}), from now on we assume the condition $|a| <1$ which corresponds to the ellipticity of the top order coefficients and also $\alpha<2$ as explained in the Introduction.

We start by the $L^2$ theory. We use the Sobolev spaces of Section \ref{Weighted}  and also  $H^{1}_c:=\{u \in L^2_{c} : \nabla u \in L^2_{c}\}$ equipped with the inner product
\begin{align*}
	\left\langle u, v\right\rangle_{H^1_{c}}:= \left\langle u, v\right\rangle_{L^2_{c}}+\left\langle \nabla u, \nabla v\right\rangle_{L^2_{c}}.
\end{align*}
and  consider the  form  in $L^2_{c-\alpha}$ with $D(\mathfrak{a})=H_c^1\cap L^2_{c-\alpha}\subset L^2_{c-\alpha}$ and 
\begin{align*}
	\mathfrak{a}(u,v)
	:&=
	\int_{\R_+}  D_y u\, D_y \overline{v}\,y^{c} \,dy-2ia\cdot\xi\int_{\R_+} (D_yu)\, \overline{v}\,y^{c} \,dy+\int_{\R_+} |\xi|^2u \overline{v} y^{c} \,dy , \\[1ex]
	&=\int_{\R_+}  y^\alpha D_y u  D_y \overline{v}\,y^{c-\alpha} \,dy-2ia\cdot\xi\int_{\R_+} y^{\alpha}D_yu\, \overline{v}\,y^{c-\alpha} \,dy +\int_{\R_+} |\xi|^2 y^{\alpha} u\overline{v} \,y^{c-\alpha} \,dy. 
\end{align*}
The more pedantic second line above writes the form with respect to the  reference measure $y^{c-\alpha} dy$, rather than $y^c dy$.

We define $L$ in $L^2_{c-\alpha}$ as the operator  associated to the form $\mathfrak{a}$, that is
\begin{align*} 
\nonumber D(L)&=\{u \in D(\mathfrak{a}): \exists  f \in L^2_{c-\alpha} \ {\rm such\ that}\  \mathfrak{a}(u,v)=\int_0^\infty f \overline{v}y^{c-\alpha}\, dy\ {\rm for\ every}\ v\in D(\mathfrak{a})\},\\[1ex]
  L u&=-f.
\end{align*}
If $u,v$ are smooth functions with compact support, it is easy to see integrating by parts that $$-\mathfrak a (u,v)= \langle y^\alpha(B u+2ia\cdot \xi D_y u-|\xi|^2u), \overline v\rangle_{L^2_{c-\alpha}}, $$
so that 
$ L$ is a realization of $y^\alpha(B +2ia\cdot \xi D_y-|\xi|^2)$. 

\medskip
In the next lemmas we use two isometries which transform the operator into a simpler form. 

The first one is $Tu(y)=\left(1-\frac{\alpha}{2}\right)^\frac{1}{2}u(y^{1-\frac \alpha 2})$ (this is the map $T_{-\frac{\alpha}{2}}$ of  Appendix B) and allows to remove the power $y^\alpha$ in front of the Bessel operator $B$ thus getting an equivalent operator $\left(\beta+1\right)^{-2}	\tilde L$ defined by 
 \begin{align*}
	\tilde L=\tilde B+ib(\beta+1) y^{\beta}D_y-(\beta+1)^2|\xi|^2 y^{2\beta},\qquad \beta=\frac{\alpha}{2-\alpha}
 \end{align*}

Here
\begin{align*}
\tilde B=D_{yy}+\frac{\tilde c}{y} D_y,\qquad \tilde c=\frac{2c-\alpha}{2-\alpha},\qquad    b=2a\cdot\xi
\end{align*}
and, by the assumption $\alpha<2$,  one has $\beta+1>0$ and moreover  $\tilde c+1>0$ if and only if $c+1-\alpha>0$.

This is contained  in Proposition \ref{Isometry action}  at level of operators but we state it below in the language of forms, omitting the elementary computations. 

\begin{lem}\label{Lemma forme eq 1}
Let  $c+1-\alpha>0$. Setting  $\tilde c=\frac{2c-\alpha}{2-\alpha}>-1$,   let us consider the isometry
$$T=:L^2_{\tilde c}\to L^2_{c-\alpha},\quad  Tu(y)=\left(1-\frac{\alpha}{2}\right)^\frac{1}{2}u(y^{1-\frac \alpha 2}).$$
 Then, with  $b=2a\cdot\xi$, $\beta=\frac{\alpha}{2-\alpha}$, one has   $T\left(H_{\tilde c}^1\cap L^2_{\tilde c+2\beta}\right)=H^1_c \cap L^2_{c-\alpha}$ and  
$$
\mathfrak{a}(u,v)=	\left(\beta+1\right)^{-2}\mathfrak{\tilde a}(T^{-1}u,T^{-1}v)
,\qquad u, v\in H^1_c \cap L^2_{c-\alpha}$$
where $D(\mathfrak{\tilde a})=H_{\tilde c}^1\cap L^2_{\tilde c+2\beta }$
\begin{align*}
	\mathfrak{\tilde a}(u,v)=
\int_{\R_+}  D_y u\, D_y \overline{v}\,y^{\tilde c} \,dy-ib(\beta+1)\int_{\R_+} y^\beta (D_yu)\, \overline{v}\,y^{\tilde c} \,dy+|\xi|^2 (\beta+1)^2\int_{\R_+} y^{2\beta}u \overline{v} y^{\tilde c} \,dy.
\end{align*}
\end{lem}
We introduce now the quadratic form 
\begin{align}\label{est Qa}
	Q_{a}(\xi)=|\xi|^2-|a\cdot\xi|^2,\qquad(1-|a|^2)|\xi|^2\leq Q_{a}(\xi)\leq |\xi|^2,\qquad  \xi\in\R^N.
\end{align}
A second isometry $S$ removes  the term $ib(\beta+1) y^{\beta}D_y$ from the operator $\tilde L$ introducing a complex potential. This leads to the operator $\left(\beta+1\right)^{-2}\tilde A_{\tilde b, \beta}$ defined by 
\begin{equation} \label{Btilde}
\tilde A_{\tilde b, \beta}=\tilde B-i\frac{ \tilde b(\tilde c+\beta)}{2} y^{\beta-1}-(\beta+1)^2Q_a(\xi)y^{2\beta}
\end{equation}
where $\tilde c= \frac{c-\alpha}{2-\alpha}$,  $\tilde b=b(\beta+1)=2a\cdot\xi (\beta+1)$.

\begin{lem} \label{iso}
With the notation above let  us consider the isometry
\begin{align*}
	S:L^2_{\tilde c}\to L^2_{\tilde c}, \qquad (Su)(y)=e^{-i\frac{ \tilde b}{2(\beta +1)} y^{\beta+1}}u(y)=e^{-ia\cdot\xi y^{\frac{2}{2-\alpha}}}u(y).
\end{align*}
If $D(\mathfrak{\tilde a})=H_{\tilde c}^1\cap L^2_{\tilde c+2\beta}$,  then  one has $S( D(\mathfrak{\tilde a}))= D(\mathfrak{\tilde a})$ and 
\begin{align*}
\mathfrak{\tilde a}(u, v)=\mathfrak{\tilde a}_{\tilde b}(S^{-1}u, S^{-1}v),\qquad u,\ v\in D(\mathfrak{\tilde a}).
\end{align*}
\begin{align*}
	\nonumber \mathfrak{\tilde a}_{\tilde b}(u,v)&=
	\int_{\R_+}
	D_y u D_y \overline{v}y^{\tilde c}\,dy
	- i\frac  {\tilde b} 2\int_{\R_+} D_y\left(u\overline{v}\right)y^{\tilde c+\beta}\,dy+(\beta+1)^2Q_a(\xi)\int_{\R_+} u\overline{v}y^{\tilde c+2\beta}\,dy.
	\end{align*}
\end{lem}


The above lemmas say that the operator $y^\alpha L_{2a\cdot\xi}-y^\alpha|\xi|^2$ and the associated form $\mathfrak{a}$ are equivalent, by mean of the isometry $T\circ S$, to $(\beta+1)^{-2}\tilde A_{ \tilde b, \beta}$ and $(\beta+1)^{-2}\mathfrak{\tilde a}_{\tilde b}$ and motivates the next section.

\begin{os} \label{parametri} {\rm Let us explain briefly the restrictions on the parameters $c, \alpha, \beta$ which appear in this section. The operator $B=D_{yy}+\frac cy D_y$ is considered here only for $c>-1$, so that the measure $y^c dy $ is finite in a neighborhood of $0$, to impose Neumann boundary conditions. The case $c \leq -1$ can be also considered under Dirichlet boundary conditions, see \cite{MNS-Caffarelli}. Note also that the set $\mathcal D$ defined in (\ref{defD}) is contained in the domain of the form associated with $B$, if and only if $c>-1$. All other restrictions on the parameters come from this choice. For example, the condition $c+1-\alpha>0$  in Lemma \ref{Lemma forme eq 1} is equivalent to $\tilde c+1>0$. 
} 

\end{os}

\subsection{The auxiliary operator $A_{b,\beta}=B-i\frac{ b( c+\beta)}{2} y^{\beta-1}-(\beta+1)^2Q_a(\xi)y^{2\beta}.$}\label{section A_b}
As explained in the above remark, we always assume that 
\begin{align*}
|a|<1, \quad c+1>0, \quad 	\beta+1>0.
\end{align*}
Setting $b=2a\cdot \xi(\beta+1),\ Q_a(\xi)=|\xi|^2-|a\cdot \xi|^2$, we 
consider the form  $\mathfrak{\tilde a}_b$  defined on $$D(\mathfrak{\tilde a}_b)=H^1_c\cap L^2_{c+2\beta}\subset L^2_c$$ by
\begin{align}\label{form tilde}
	\nonumber &\mathfrak{\tilde a}_b(u,v)=\langle D_y u, D_y v \rangle_{L^2_c}- i\frac  b 2\left\langle u,D_yv\right\rangle_{L^2_{c+\beta}}- i\frac  b 2\left\langle D_yu,v\right\rangle_{L^2_{c+\beta}}+(\beta+1)^2Q_a(\xi)\left\langle u,v \right\rangle_{L^2_{c+2\beta}}\\[1.5ex]&=
	\int_{\R_+}
	D_y u D_y \overline{v}y^c\,dy
	- i\frac  b 2\int_{\R_+} D_y\left(u\overline{v}\right)y^{c+\beta}\,dy+(\beta+1)^2Q_a(\xi)\int_{\R_+} u\overline{v}y^{c+2\beta}\,dy
\end{align}
and its   associated operator $A_{b,\beta}$ in $L^2_c$. 
Since for smooth functions with compact support away from the origin
$$
(c+\beta)\int_0^\infty u \overline{v} y^{c+\beta -1}\, dy=
\int_0^\infty u\big  (D_y (y^{c+\beta}\overline{v})-y^{c+\beta}  D_y \overline{v})\big )\, dy=- \int_0^\infty D_y(u\overline{v})y^{c+\beta}\, dy,
$$
the operator $A_{b,\beta}$ is defined on smooth functions by 
\begin{align*}
	A_{b,\beta}:= B-i\frac{b(c+\beta)}{2}y^{\beta-1}-(\beta+1)^2Q_a(\xi)y^{2\beta}, \quad  A_0=B-(\beta+1)^2Q_a(\xi)y^{2\beta}.
\end{align*}

We collect in the following proposition the main properties satisfied by $\mathfrak{\tilde a}_b$.

\begin{prop}\label{prop form a_b tilde}
	The form  $\mathfrak{\tilde a}_b$ is  accretive  and closed in $L^2_{c}$. Moreover
	\begin{itemize}
		\item[(i)] the adjoint form $\mathfrak{\tilde a}_{b}^\ast:(u,v)\mapsto \overline{\mathfrak{\tilde a}_b(v,u)}$ satisfies $\mathfrak{\tilde a}_{b}^\ast =\mathfrak{\tilde a}_{-b}$; 
		\item[(ii)] its real part  is the positive form
		 $$\Rp{\mathfrak{\tilde a}_b}(u,v):=\frac{\mathfrak{\tilde a}_b(u,v)+\mathfrak{\tilde a}_b^\ast (u,v)}{2}=\langle D_y u,D_y v\rangle_{L^2_c}+(\beta+1)^2Q_a(\xi)\left\langle u ,  v  \right\rangle_{L^2_{c+2\beta}};$$ 
		\item[(iii)] for any $u \in H^1_c\cap L^2_{c+2\beta}$
	\begin{align*}
		\left|\Ip{\mathfrak{\tilde a}_b}(u,u)\right|&\leq \frac{|b|}{(\beta+1)Q_a(\xi)^\frac{1}{2}} \Rp{\mathfrak{\tilde a}_b}(u,u)=\frac{|a|}{\sqrt{1-|a|^2}} \Rp{\mathfrak{\tilde a}_b}(u,u).
	\end{align*}
	\end{itemize}
\end{prop}
{\sc Proof.}  Properties (i) and (ii) are immediate consequences of the definition.	Since $\Rp{\mathfrak{\tilde a}_b(u,u)}=\|D_y u\|^2_{L^2_c}+(\beta+1)^2Q_a(\xi)\|u\|_{L^2_{c+2\beta}}^2\geq 0$,    $\mathfrak{\tilde a}_b$ 
is  accretive and, furthermore,  the norm induced by the form $\mathfrak{\tilde a}_b$  coincides with the one of  $H^1_c\cap L^2_{c+2\beta}$ and then $\mathfrak{\tilde a}_b$ is closed.  

To prove   (iii), we  use Young's inequality and the elementary identity $D_y(|u|^2)=2\Rp\left(\overline u D_yu\right)$ for  $u\in H^1_c$. Then   
\begin{align*}
	\left|\Ip{\mathfrak{\tilde a}_b}(u,u)\right|&=\left|-\frac{b}2\int_0^\infty D_y\left(|u|^2\right)y^{c+\beta}\,dy\right|=\left|-b\int_0^\infty \Rp\left(\overline u D_yu\right)y^{c+\beta}\,dy\right|\\[1.5ex]
	&\leq |b| \left(\|D_yu\|_{L^2_c} \|u\|_{L^2_{c+2\beta}}\right)=\frac{ |b|}{(\beta+1)Q_a(\xi)^\frac{1}{2}} \left(\|D_yu\|_{L^2_c}\; {(\beta+1)Q_a(\xi)^\frac{1}{2}}\|u\|_{L^2_{c+2\beta}}\right) \\[1.5ex]
	&\leq  \frac{ |b|}{2(\beta+1)Q_a(\xi)^\frac{1}{2}}\left(\|D_yu\|^2_{L^2_c}+(\beta+1)^2Q_a(\xi)\|u\|^2_{L^2_{c+2\beta}}\right)\\[1.5ex]
	&= \frac{|b|}{2(\beta+1)Q_a(\xi)^\frac{1}{2}}  \Rp{\mathfrak{\tilde a}_b(u,u)}=\frac{|a\cdot \xi|}{Q_a(\xi)^\frac{1}{2}}\Rp{\mathfrak{\tilde a}_b(u,u)}\leq \frac{|a|}{\sqrt{1-|a|^2}}\Rp{\mathfrak{\tilde a}_b(u,u)}.
\end{align*}

By standard theory on sesquilinear forms we have  the following results.

\begin{prop}\label{kernel Vreal tilde1}
	The   operator $A_{b,\beta}$ generates  an  analytic semigroup  of angle $\frac\pi 2-\arctan \frac{|a|}{\sqrt{1-|a|^2}}$ in $L^2_c$   which satisfies
	\begin{align*}
		\|e^{zA_{b,\beta}}f\|_{L^2_c}\leq \|f\|_{L^2_c},\qquad \forall z\in \Sigma_{\frac\pi 2-\arctan \frac{|a|}{\sqrt{1-|a|^2}} }.
	\end{align*} 
	Moreover
	\begin{itemize}
		\item[(i)] The semigroup $\left(e^{tA_{b,\beta}}\right)_{t\geq 0}$ is  $L^\infty$-contractive  and it is dominated by $e^{tB}$, that is 
		\begin{align*}
			|e^{tA_{b,\beta}}f|\leq e^{tB}|f|,\qquad t>0,\quad  f\in L^2_c.
		\end{align*}
		\item[(ii)] $\left(e^{tA_{b,\beta}}\right)_{t\geq 0}$ is a semigroup of integral operators and its heat kernel $\tilde p_{b,\beta}$, taken with respect to the measure $\rho^cd\rho$, satisfies for some constant $C$ independent of $b,\ \beta$
		\begin{align*}
			|\tilde p_{b,\beta}(t,y,\rho)|\leq C t^{-\frac{1}{2}}  \rho^{-c}\left (\frac{\rho}{t^{\frac{1}{2}}}\wedge 1 \right)^{c}
			\exp\left(-\frac{|y-\rho|^2}{\kappa t}\right),\qquad \text{for a.e. $y,\rho>0$}.
		\end{align*}     
		\item[(iii)] ${A}^\ast_{b,\beta}= A _{-b,\beta}$.
		and for any $s>0$ the operator satisfies the scaling property
		\begin{align*}
			I_{\frac 1 s}\circ A_{b,\beta} \circ I_s= s^2\left(B-i\frac{ b( c+\beta)}{2s^{1+\beta}} y^{\beta-1}-\frac{(\beta+1)^2}{s^{2\beta+2}}Q_a(\xi)y^{2\beta}\right),\qquad I_s u(y):=u(sy). 
		\end{align*}
	\end{itemize}  
\end{prop}
{\sc Proof.} The generation properties  follows  using  Proposition \ref{prop form a_b tilde} and \cite[Teorems 1.52, 1.53]{Ouhabaz}. 	To prove  (i) we observe, preliminarily,  that the operator $B$ is associated with the form $\mathfrak b(u,v)=\langle D_y u, D_y v \rangle_{L^2_c}$ and its generated semigroup $e^{tB}$ is sub-Markovian  since $\mathfrak b$ satisfies the hypotheses of  \cite[Corollary 2.17]{Ouhabaz}.  The domination property for $e^{tA_{b,\beta}}$ then  follows  from \cite[Theorem 2.21]{Ouhabaz}. In particular  $e^{tA_{b,\beta}}$ inherits the  $L^\infty$-contractivity of $e^{tB}$. (ii) is a consequence of  \cite[Proposition 1.9]{ArendtBukhalov} since 
$e^{tA_{b,\beta}}$ is dominated by the positive integral operator $e^{tB}$ whose kernel satisfies the stated estimate, see  \cite[Proposition 2.8]{MNS-Caffarelli} where, however, the  kernel is written with respect to the Lebesgue measure. (iii) follows from (i) of Proposition \ref{prop form a_b tilde} and by elementary computations.\qed

As in \cite[Section 5]{MNS-Singular-Half-Space}, we can  extend the above heat kernel estimates to complex times. 

\begin{teo}\label{kernel Vreal tilde}
For every $0\leq\nu< \frac\pi 2-\arctan \frac{|a|}{\sqrt{1-|a|^2}}$ the   heat kernel $\tilde p_{b,\beta}$, taken with respect to the measure $\rho^cd\rho$, satisfies for some constant $C_\nu$ independent of $b,\ \beta$
		\begin{align*}
			|\tilde p_{b,\beta}(z,y,\rho)|\leq C_\nu |z|^{-\frac{1}{2}}  \rho^{-c}\left (\frac{\rho}{|z|^{\frac{1}{2}}}\wedge 1 \right)^{c}
			\exp\left(-\frac{|y-\rho|^2}{\kappa |z|}\right),
			\end{align*}     
		for a.e. $y,\rho>0$, $\forall z\in \Sigma_{\nu}$.
\end{teo}

\subsection{Generation properties and domain characterization}
Generation properties and  kernel estimates  for the original operator $$y^\alpha L_{2a\cdot\xi}-y^\alpha|\xi|^2=y^\alpha (B+i2a\cdot\xi D_y-|\xi|^2)$$ can be deduced
 by the analogous properties of the auxiliary operator $A_{b,\beta}$ of Section \ref{section A_b}. Indeed from Lemmas  \ref{Lemma forme eq 1}, \ref{iso} we have
\begin{align*}
	\mathfrak{a}(u,v)=(\beta+1)^{-2}\,\mathfrak{\tilde a}_{\tilde b}(S^{-1}T^{-1}u,S^{-1}T^{-1}v ),\qquad u,v\in D(\mathfrak a).
\end{align*}

This implies that
\begin{align}\label{Equivalence L_b-A_b}
y^\alpha L_{2a\cdot\xi}-y^\alpha|\xi|^2=(T\circ S)\Big[(\beta+1)^{-2}\tilde A_{\tilde b, \beta}\Big](T\circ S)^{-1},
\end{align}
where $\beta=\frac {\alpha}{2-\alpha}$, $ \tilde b=2a\cdot \xi(\beta+1)$, $\tilde c=\frac{2c-\alpha}{2-\alpha}$ and
\begin{align*} 
\tilde A_{\tilde b, \beta}=\tilde B-i\frac{ \tilde b(\tilde c+\beta)}{2} y^{\beta-1}-(\beta+1)^2Q_a(\xi)y^{2\beta}.
\end{align*}
Note that, by construction and by Proposition \ref{Isometry action} we have
\begin{align}\label{Equivalence L_b-A_b 2}
	\tilde B=D_{yy}+\frac{\tilde c}{y} D_y,\qquad y^\alpha B=T\Big[(\beta+1)^{-2}\tilde B\Big] T^{-1}.
\end{align}
\begin{teo} \label{kernel V}
Let  $c+1-\alpha>0$. Then the operator $y^\alpha L_{2a\cdot\xi}-y^\alpha|\xi|^2$ generates a contractive analytic semigroup  of angle $\frac\pi 2-\arctan \frac{|a|}{\sqrt{1-|a|^2}}$ in $L^2_{c-\alpha}$
Moreover
\begin{itemize}
	\item[(i)] The semigroup $\left(e^{ty^\alpha\left( L_{2a\cdot\xi}-|\xi|^2\right)}\right)_{t\geq 0}$  is dominated by $e^{ty^\alpha B}$, that is 
	\begin{align*}
		|e^{ty^\alpha\left( L_{2a\cdot\xi}-|\xi|^2\right)}f|\leq e^{ty^\alpha B}|f|,\qquad t>0,\quad  f\in L^2_{c-\alpha}.
	\end{align*}
	\item[(ii)] $\left(e^{ty^\alpha\left( L_{2a\cdot\xi}-|\xi|^2\right)}\right)_{t\geq 0}$ is a semigroup of integral operators and its heat kernel  $p_\alpha$, taken with respect to the measure $\rho^{c-\alpha}d\rho$, satisfies for some constant $C_\nu$
	$$
	|p_\alpha(z,y,\rho)|\leq C_\nu|z|^{-\frac {1} 2}\, \rho^{-\frac\alpha 2} \left (\frac{\rho}{|z|^{\frac{1}{2-\alpha}}}\wedge 1 \right)^{c+\frac\alpha 2}
	\exp\left(-\frac{|y^{1-\frac\alpha 2}-\rho^{1-\frac\alpha 2}|^2}{\kappa |z|}\right),
	$$
	for a.e. $y,\rho>0$, $\forall z\in \Sigma_{\nu}$, $0\leq \nu<\frac\pi 2-\arctan \frac{|a|}{\sqrt{1-|a|^2}}$.
\end{itemize}  
\end{teo}
{\sc{Proof.}}
The proof  is  simply a translation of the results for $	\tilde A_{\tilde b, \beta}$ of  Proposition \ref{kernel Vreal tilde1} and of Theorem \ref{kernel Vreal tilde} by using the identity \eqref{Equivalence L_b-A_b}. For example,  (i) follows since, by construction, we have for any $g $, $|TSg|=T|g|$ and therefore using  (i) of Proposition  \ref{kernel Vreal tilde1} and \eqref{Equivalence L_b-A_b}, \eqref{Equivalence L_b-A_b 2} we get  for $t>0$, $f\in L^2_{c-\alpha}$
\begin{align*}
	|e^{ty^\alpha\left( L_{2a\cdot\xi}-|\xi|^2\right)}f|&=|TSe^{t(\beta+1)^{-2}	\tilde A_{\tilde b, \beta}}S^{-1}T^{-1}f|=T|e^{t(\beta+1)^{-2}\tilde A_{\tilde b, \beta}}S^{-1}T^{-1}f|\\[1ex]
	&\leq Te^{t(\beta+1)^{-2}\tilde B}|S^{-1}T^{-1}f|=Te^{t(\beta+1)^{-2}\tilde B}T^{-1}|f|=e^{ty^\alpha B}|f|.
	\end{align*}
  \qed

Now we prove that the semigroup $e^{z (y^\alpha(L_{2a\cdot\xi}-|\xi|^2))}$ extrapolates to the spaces $L^p_m$.

\begin{prop}\label{Gen DN}
	If  $1<p < \infty$ and  $0< \frac{m+1}{p} <c+1-\alpha$, then  $(e^{z(y^\alpha(L_{2a\cdot\xi}-|\xi|^2))})$ is an  analytic semigroup of angle $\frac\pi 2-\arctan \frac{|a|}{\sqrt{1-|a|^2}}$ in $L^p_m$. 
\end{prop}
{\sc Proof. } All properties for  $p=2$, $m=c-\alpha$ are contained in Theorem  \ref{kernel V}.  The boundedness of $e^{z(y^\alpha(L_{2a\cdot\xi}-|\xi|^2))}$ in $L^p_m$ then  follows from  \cite[Proposition 12.2]{MNS-Caffarelli}.
The semigroup law is inherited  from the one of $L^2_{c-\alpha}$ via a density argument and we have only to prove the  strong continuity at $0$.  Let $f, g \in C_c^\infty (\R_+)$. Then as $z \to 0$, $z \in \Sigma_{\frac\pi 2-\arctan \frac{|a|}{\sqrt{1-|a|^2}}}$, 
\begin{align*}
&\int_0^\infty (e^{z (y^\alpha(L_{2a\cdot\xi}-|\xi|^2))}f)\, g\, y^m dy=\int_0^\infty (e^{z (y^\alpha(L_{2a\cdot\xi}-|\xi|^2))}f) \,g\, y^{m-c+\alpha}y^{c-\alpha}  dy\\& \to \int_0^\infty fgy^{m-c+\alpha}y^{c-\alpha}dy  =\int_0^\infty fgy^{m}dy,
\end{align*}
by the strong continuity of $e^{z (y^\alpha(L_{2a\cdot\xi}-|\xi|^2))}$ in $L^2_{c-\alpha}$. Let us observe now that, using  Theorem \ref{kernel V} and  \cite[Proposition 12.2]{MNS-Caffarelli},   the family $\left\{e^{z (y^\alpha(L_{2a\cdot\xi}-|\xi|^2))}, z\in \Sigma_{\frac\pi 2-\arctan \frac{|a|}{\sqrt{1-|a|^2}}}\right\}$  is uniformly bounded on $\mathcal B(L^p_m)$. By   density,  the previous limit holds for every $f \in L^p_m$, $g \in L^{p'}_m$. The semigroup is then weakly continuous, hence strongly continuous.
\qed

Following the same lines of  \cite[Theorem 7.2]{MNS-PerturbedBessel}, we get the  following  $\mathcal R$-boundedness result (see Appendix A for the relevant definitions). 
\begin{cor}\label{Rbound lR_l}
 Let $1<p<\infty$ such that $0<\frac{m+1}p<c+1-\alpha$. Then  the following properties hold.
 For every $ 0\leq \nu< \pi -\arctan \frac{|a|}{\sqrt{1-|a|^2}}$  the families of operators  
\begin{align*}
\left\{e^{zy^\alpha(L_{2a\cdot\xi}-|\xi|^2)}:\; \xi\in\R^N\setminus\{0\},\ z\in\Sigma_\nu\right\},\\[1ex]
\left\{\lambda\left(\lambda-y^\alpha(L_{2a\cdot\xi}-|\xi|^2) \right)^{-1}:\;\xi\in\R^N\setminus\{0\},\ \lambda \in \Sigma_\nu\right\}
\end{align*}
 are $\mathcal R$-bounded in $L^p_m$.
 \end{cor}
{\sc{Proof.}} The proof follows almost identically to that of   \cite[Theorem 7.2]{MNS-PerturbedBessel} since from (ii) of Theorem \ref{kernel V} one has (using the notation of  \cite[Theorem 7.2]{MNS-PerturbedBessel}) for every $\leq \nu<\frac\pi 2-\arctan \frac{|a|}{\sqrt{1-|a|^2}}$  and for some positive constant $C$ 
\begin{align*}
	\left|e^{zy^\alpha(L_{2a\cdot\xi}-|\xi|^2)}f\right|\leq CS_\alpha^{-c}(|z|)|f|,\quad f\in L^p_{m},\quad  z\in\Sigma_\nu.
\end{align*}
The $\mathcal R$-boundedness on $\xi\in\R^N\setminus\{0\}$ follows since  the right hand side  does not depend on $\xi$.
\qed 

In our investigations of degenerate  $N$-d problems, we need,  in the case $\alpha=0$, to add a potential having non-negative real part to the operator of the latter proposition; this force to deal only with the  semigroup on the real axis.

Let $V\in L^1_{loc}\left(\R^+,y^{c}\, dy\right)$ be a   potential having non-negative real part and let   $L_{2a\cdot\xi}-V$  be the operator in $L^2_c$ associated with the form
\begin{align*}
	&\mathfrak{a}_{V }(u,v)=
	\int_{\R_+}
	\left( D_y u D_y \overline{v}-2ia\cdot\xi D_y u  \overline{v}+Vu\overline{v}\right)y^{c}\,dy
\end{align*}
defined on the  domain
\begin{align*}
	\mathcal{F}:=H^1_c\cap L^2\left(\R_+,Vy^cdy\right).
\end{align*}

\begin{prop}\label{kernel Vreal}
Let $V\in L^1_{loc}\left(\R^+,y^{c}\, dy\right)$ be a   potential having non-negative real part. Then for any $1<p<\infty$ such that $0<\frac{m+1}p<c+1$,  $L_{2a\cdot\xi}-|\xi|^2-V$ generates a $C_0$-semigroup on $L^p_m$. The generated semigroup consists of integral operators and the following estimate holds
     \begin{align*}
   \left|e^{t(L_{2a\cdot\xi}-|\xi|^2-V)}f\right|&\leq e^{t (L_{2a\cdot\xi}-|\xi|^2)}|f|\leq e^{t B}|f|,\hspace{13ex} f\in L^p_m,\quad t\geq 0.
   \end{align*}
Moreover  the families of operators
\begin{align}\label{R bound Resolv V}
	\nonumber
	\left\{e^{t(L_{2a\cdot\xi}-|\xi|^2-V)}:\; t\geq 0,\; \xi\in\R^N\setminus\{0\},\ |a|<1,\  V\in L^1_{loc}\left(\R^+,y^{c}\right), \ \Rp V\geq 0\right\},\\[1ex]
	\left\{\lambda\left(\lambda-L_{2a\cdot\xi}+|\xi|^2+V\right)^{-1}:\; \lambda >0,\;  \xi\in\R^N\setminus\{0\},\ |a|<1,\ V\in L^1_{loc}\left(\R^+,y^{c}\right), \ \Rp V\geq 0\right\}
\end{align}
are $\mathcal R$-bounded in $L^p_m$. In particular
\begin{align*}	\left\{|\xi|^2\left(|\xi|^2-L_{2a\cdot\xi}+V\right)^{-1}:\;   \xi\in\R^N\setminus\{0\},\ |a|<1,\ V\in L^1_{loc}\left(\R^+,y^{c}\right), \ \Rp V\geq 0\right\}
\end{align*} 
is $\mathcal R$-bounded in $L^p_m$.
\end{prop}
{\sc{Proof.}}
The generation results can be proved as  in  Proposition \ref{kernel Vreal tilde1}. If $\mathfrak a$ is the  form associated with $L_{2a\cdot\xi} -|\xi|^2$, then  $L_{2a\cdot\xi}-|\xi|^2-V$ is associated to ${\mathfrak a}_V:=\mathfrak a(u,v)+\langle Vu,v\rangle_{L^2_c}$ and,  by the standard theory on sesquilinear forms,  $L_{2a\cdot\xi}-V$  generates a $C_0$-semigroup on $L^2_c$. 
The domination properties follow from  \cite[Theorem 2.21]{Ouhabaz} again.  The extrapolation on $L^p_m$ follows as in Proposition \ref{Gen DN}. The $\mathcal R$-boundedness of the semigroup  follows by domination from the $\mathcal R$-boundedness of $e^{t (L_{2a\cdot\xi}-|\xi|^2)}$ using Corollary \ref{Rbound lR_l} with $\alpha=0$. The $\mathcal R$-boundedness of the resolvent family follows from   Proposition \ref{Mean R-bound} by writing the resolvent as the Laplace transform of the semigroup. The last claim follows by simply   specializing \eqref{R bound Resolv V} taking $\lambda=|\xi|^2$.
\qed

We now prove that  the domain of $y^\alpha L_{2a\xi}+|\xi|^2y^\alpha$ is $W^{2,p}_{\mathcal N}(\alpha,m)\cap  L^p_{m+\alpha p}$, under slightly more restrictive hypotheses than those of Proposition \ref{Gen DN}. Indeed, in what follows, we assume, besides $c+1-\alpha>0$ also the condition $c+1>0$.  
\begin{lem}\label{Lemma Rb 3}
	Assume that   $c+1>0$ and  $c+1-\alpha>0$. If  $ \lambda \in \mathbb C^+$ and $\xi\in \R^N\setminus\{0\}$,  then 
	\begin{align*}
	\left(\lambda-y^\alpha L_{2a\xi}+|\xi|^2y^\alpha\right)^{-1}f=\left(|\xi|^2-L_{2a\cdot \xi} +\frac{\lambda}{y^\alpha}\right)^{-1}\left(\frac{f}{y^\alpha}\right),\qquad \forall f\in C^\infty_c((0,\infty)).
\end{align*}	
\end{lem} 
{\sc{Proof.}} Under the assumptions  $y^\alpha L_{2a\xi}-|\xi|^2 y^\alpha$ and $L_{2a\xi}-\lambda y^{-\alpha}$  generate a semigroup  on $L^2_{c-\alpha}$ and $L^2_c$, respectively, see Theorem \ref{kernel V} and Proposition \ref{kernel Vreal}. Since $\Rp \lambda>0$,  both resolvents are well defined but act in different spaces.

Let  $\mathfrak{a}$ , $\mathfrak{a}_{\lambda y^{-\alpha}}$ be the forms  associated   to $y^\alpha L_{2a\xi}-|\xi|^2 y^\alpha$   in $L^2_{c-\alpha}$ and  $L_{2a\xi}-\lambda y^{-\alpha}$ in  $L^2_c$
\begin{align*}
	&\mathfrak{a}(u,v)=
	\int_{\R_+}
	\left( D_y u D_y \overline{v}+2ia\cdot\xi D_y u  \overline{v}+|\xi|^2 u\overline{v}\right)y^{c}\,dy,\\&\mathfrak{a}_{ \lambda y^{-\alpha} }(u,v)=
	\int_{\R_+}
	\left( D_y u D_y \overline{v}+2ia\cdot\xi D_y u  \overline{v}+\lambda y^{-\alpha}u\overline{v}\right)y^{c}\,dy.
\end{align*}
They are defined on the common domain
\begin{align*}
	\mathcal{F}:=\left\{u\in L^2_{c-\alpha}\cap L^2_c:D_y u\in L^2_{c}\right\}
\end{align*}

 Given   $f\in C^\infty_c((0,\infty))$ let  $u:= \left(|\xi|^2-L_{2a\cdot\xi} +\frac{\lambda}{y^\alpha}\right)^{-1}\left(\frac{f}{y^\alpha}\right)$. In order to prove  that the equality  $u=\left(\lambda-y^\alpha L_{2a\cdot\xi} +|\xi|^2 y^\alpha\right)^{-1}f$ holds,  we have to show that $u\in \mathcal{F}$ and that for every $v\in \mathcal F$, $u$ satisfies the weak equality 
 \begin{align}\label{Rb 3 eq1}
&\int_0^\infty f \overline{v}y^{c-\alpha}\, dy=\int_0^\infty \lambda u \overline{v}y^{c-\alpha}\, dy+\mathfrak{a}_{\alpha,|\xi|^2 y^\alpha}(u,v)\\&=\int_0^\infty (\lambda y^{-\alpha}u \overline{v}+ D_y u D_y \overline{v}+2ia\cdot\xi D_y u  \overline{v}+|\xi|^2  u\overline{v})y^{c}\, dy.
 \end{align}
 By construction  $u$ is in the domain  of $L_{2a\cdot\xi}-\lambda y^{-\alpha}$ which is contained in $\mathcal F$ and satisfies  
 \begin{align*}
	&\int_0^\infty \frac{f}{y^\alpha} \overline{v}y^{c}\, dy=\int_0^\infty |\xi|^2 u \overline{v}y^{c}\, dy+\mathfrak{a}_{\alpha,\lambda y^{-\alpha}}(u,v)\\&=\int_0^\infty (|\xi|^2 u \overline{v}+ D_y u D_y \overline{v}+2ia\cdot\xi D_y u  \overline{v}+\lambda y^{-\alpha} u\overline{v})y^{c}\, dy,
\end{align*}
which is the same as  \eqref{Rb 3 eq1}. \qed

\begin{os} \label{doublecond} {\rm In the next result we relate the resolvent of $y^\alpha L_{2a\cdot\xi}- y^\alpha$ with that of $L_{2a\cdot\xi}-\frac{1}{y^\alpha}$, in the sense specified below.  We shall assume both the conditions  $0<\frac{m+1}p<c+1-\alpha$ and $-\alpha<\frac{m+1}p<c+1-\alpha$ (that is $\alpha^-<\frac{m+1}p<c+1-\alpha$). The first guarantees that $y^\alpha L_{2a\cdot\xi}- y^\alpha$ is a generator in $L^p_{m}$ and the second that  $L_{2a\cdot\xi}-\frac 1{y^\alpha}$ is a generator in $L^p_{m+\alpha p}$.}
\end{os}

\begin{cor}\label{Cor Rb}
	Assume that  $\alpha^-<\frac{m+1}p<c+1-\alpha$. If $ \lambda \in \mathbb C^+$ and $\xi\in\R^N\setminus\{0\}$, then 
	\begin{itemize}
		\item[(i)] for every $f\in L^p_m$ 
		\begin{align*}
			\left(\lambda-y^\alpha L_{2a\cdot\xi}+|\xi|^2 y^\alpha\right)^{-1}f=\left(|\xi|^2- L_{2a\cdot\xi}+\frac{\lambda}{y^\alpha}\right)^{-1}\left(\frac{f}{y^\alpha}\right)\in L^p_{m+\alpha p}\cap L^p_m;
		\end{align*}
		
		\item[(ii)] the operator 	$ y^\alpha\left(\lambda-y^\alpha L_{2a\cdot\xi}+|\xi|^2 y^\alpha\right)^{-1}$ is  bounded in $L^p_m$;
\item[(iii)] the operator 	$ \frac{1}{y^\alpha}\left(|\xi|^2- L_{2a\cdot\xi}+\frac{\lambda}{y^\alpha}\right)^{-1}$
		is  bounded in $L^p_{m+\alpha p}$.
	\end{itemize}
	
\end{cor} 
{\sc{Proof.}} Equality (i) is proved in Lemma \ref{Lemma Rb 3} for any $f\in C^\infty_c((0,\infty))$. Since  $\left(\lambda-y^\alpha L_{2a\cdot\xi}+|\xi|^2 y^\alpha\right)^{-1}$ is bounded form $L^p_m$ into itself and $\left(|\xi|^2- L_{2a\cdot\xi}+\frac{\lambda}{y^\alpha}\right)^{-1}\left(\frac{\cdot}{y^\alpha}\right)$ is bounded from $L^p_m$ to $L^p_{m+\alpha p}$,  by density, (i) holds for every $f \in L^p_m$.
Parts (ii), (iii) are consequence of (i).\qed

To characterize the domain of $y^\alpha(L_{2a\cdot\xi}-|\xi|^2)$ we need the following lemmas.

\begin{lem}   \label{Sobolev Formula magica}
	Let   $m\in\R$, $p>1$. Then 
	\begin{align*}
		W^{2,p}_{\mathcal N}(\alpha,m)\cap L^p_{m+\alpha p}=W^{2,p}_{\mathcal N}(0,m+\alpha p)\cap L^p_m.
	\end{align*}
\end{lem}
{\sc Proof.} We observe preliminarily that 
from \cite[Lemma 3.5]{MNS-Sobolev} (which holds for $m\in\R$ and not only for $m<2$), there exist $C>0, \eps_0>0$ such that for every  $u\in W^{2,p}_{loc}(\R_+)$ one has 
\begin{align*}
\|y^{\frac{\alpha}{2}} D_yu\|_{L^p_m((1, \infty))} \leq C \left (\eps \|y^\alpha D_{yy}u\|_{L^p_m((1, \infty))} +\frac{1}{\eps} \|u\|_{L^p_m((1, \infty))} \right ).
\end{align*}
This and the elementary inequality $y^\frac{\alpha}{2}\leq y^{\alpha-1}$, \;$y\leq 1$ grant that the term $y^{\frac{\alpha}{2}} D_yu$ can be discarded from the definition of the Sobolev space showing that
\begin{align*}
	W^{2,p}_{\mathcal N}(\alpha,m)=\left\{u\in W^{2,p}_{loc}(\R_+):\ u,\    y^{\alpha}D_{yy}u,\ y^{\alpha-1}D_{y}u\in L^p_m\right\}.
\end{align*}
In view of the latter  equality, the required identity becomes trivial. 
\qed

\begin{lem}   \label{interpolazione}
Let  $1<p < \infty$, $\alpha^{-}<\frac{m+1}{p}<c+1-\alpha$. Then there exists $C>0$ such that for every  $u\in W^{2,p}_{\mathcal N}(\alpha,m)\cap  L^p_{m+\alpha p}$, 
$$\|y^\alpha D_{y} u\|_{L^p_m}\leq C\| y^\alpha Bu\|_{L^p_m}^\frac{1}{2}\| y^\alpha u\|_{L^p_m}^\frac{1}{2}.$$ It follows that for every $\eps>0$, $\xi\in\R^N$,
$$\||\xi|y^\alpha D_{y} u\|_{L^p_m}\leq \eps\| y^\alpha Bu\|_{L^p_m}+\frac{C}{\eps}\| |\xi|^2y^\alpha u\|_{L^p_m}.$$
\end{lem}
{\sc Proof.} We apply   \cite[Lemma 5.15]{MNS-Singular-Half-Space} with $m+\alpha p$ in place of $m$ thus obtaining
\begin{align*}
	\| D_{y} u\|_{L^p_{m+\alpha p}}\leq C\| Bu\|_{L^p_{m+\alpha p}}^\frac{1}{2}\| u\|_{L^p_{m+\alpha p}}^\frac{1}{2},\qquad u\in W^{2,p}_{\mathcal N}(0,m+\alpha p).
\end{align*}
The first inequality then follows since by Lemma \ref{Sobolev Formula magica}
 \begin{align*}
	W^{2,p}_{\mathcal N}(\alpha,m)\cap L^p_{m+\alpha p}=W^{2,p}_{\mathcal N}(0,m+\alpha p)\cap L^p_m.
\end{align*}
The second one follows by Young inequality.
\qed

\begin{lem} \label{bound}
Let   $1<p < \infty$, $\alpha^{-}<\frac{m+1}{p}<c+1-\alpha$. Then there exists $C>0$, independent of $\xi \in \R^N$, $|a|<1$, such that for every  $u\in W^{2,p}_{\mathcal N}(\alpha,m)\cap  L^p_{m+\alpha p}$, $\xi\in \R^N\setminus\{0\}$, $\Rp\lambda>0$,
$$\| |\xi|^2y^\alpha u\|_{L^p_m}\leq C\|(\lambda-y^\alpha L_{2a\cdot\xi}+|\xi|^2y^\alpha)u\|_{L^p_m}.$$
\end{lem}
{\sc Proof}. 
 Let $u\in W^{2,p}_{\mathcal N}(\alpha,m)\cap  L^p_{m+\alpha p}$ and  set $f=(\lambda-y^\alpha L_{2a\cdot\xi}+|\xi|^2y^\alpha)u$. Then, by Corollary \ref{Cor Rb} (i),
\begin{align*}
\| |\xi|^2y^\alpha u\|_{L^p_m}&=\| |\xi|^2y^\alpha (\lambda-y^\alpha L_{2a\cdot\xi}+|\xi|^2y^\alpha)^{-1}f\|_{L^p_m}=\left\| |\xi|^2y^\alpha \left(\frac{\lambda}{y^\alpha}-L_{2a\cdot\xi}+|\xi|^2\right)^{-1}\frac{f}{y^\alpha}\right\|_{L^p_m}.
\end{align*}
By Proposition \ref{kernel Vreal}            
\begin{align*}
\left\| |\xi|^2y^\alpha \left(\frac{\lambda}{y^\alpha}-L_{2a\cdot\xi}+|\xi|^2\right)^{-1}\frac{f}{y^\alpha}\right\|_{L^p_m}&=\left\| |\xi|^2 \left(\frac{\lambda}{y^\alpha}-L_{2a\cdot\xi}+|\xi|^2\right)^{-1}\frac{f}{y^\alpha}\right\|_{L^p_{m+\alpha p}}\\ 
&\leq 
C\left\| \frac{f}{y^\alpha}\right\|_{L^p_{m+\alpha p}}=C\| f\|_{L^p_{m}}
\end{align*}
for some $C$ independent of $\xi$.
\qed

\begin{teo} \label{domainLb}
Let    $1<p < \infty$, $\alpha^{-}<\frac{m+1}{p}<c+1-\alpha$ and $\xi\in \R^N\setminus\{0\}$. Then the generator of $(e^{z (y^\alpha(L_{2a\cdot\xi}-|\xi|^2))})$ is the operator $y^\alpha(L_{2a\cdot\xi}-|\xi|^2)$ with domain $W^{2,p}_{\mathcal N}(\alpha,m)\cap  L^p_{m+\alpha p}$. The set
$\mathcal {D}$ defined in \eqref{defD} 
is a core.
\end{teo}
{\sc Proof. }  We fix $0 \neq \xi$ and first  prove that the equation $ u-y^\alpha L_{2a\cdot\xi}u+|\xi|^2y^\alpha u=f$, $f\in L^p_m$, is uniquely solvable in  $W^{2,p}_{\mathcal N}(\alpha,m)\cap  L^p_{m+\alpha p}$. Let $u\in W^{2,p}_{\mathcal N}(\alpha,m)\cap  L^p_{m+\alpha p}$. By Theorem \ref{1d} (i) in the first inequality and then by  Lemma \ref{interpolazione} and Lemma \ref{bound}, there exists a positive constant $C$ such that for every $\eps>0$, $0 \leq t \leq 1$
\begin{align*}
&\|u\|_{W^{2,p}_{\mathcal N}(\alpha,m)}+\||y^\alpha u\|_{L^p_m}\leq C\|u-y^\alpha Bu\|+\|y^\alpha u\|_{L^p_m}\\[1.5ex]
&\leq  C\left(\|u-y^\alpha L_{2ta\cdot\xi}u+|\xi|^2y^\alpha u\|_{L^p_m}+|2ta\cdot\xi|\|y^\alpha D_yu\|_{L^p_m}+\|(1+|\xi|^2) y^\alpha u\|_{L^p_m}\right)
\\[1.5ex]
&\leq  C\left(\|u-y^\alpha L_{2ta\cdot\xi}u+|\xi|^2y^\alpha u\|_{L^p_m}+\eps\|y^\alpha Bu\|_{L^p_m}+\frac{1}{\eps}\|(1+|\xi|^2)y^\alpha u\|_{L^p_m}\right)
\\[1.5ex]
&\leq  C\left(\|u-y^\alpha L_{2ta\cdot\xi}u+|\xi|^2y^\alpha u\|_{L^p_m}+\eps\|u\|_{W^{2,p}_{\mathcal N}(\alpha,m)}+\frac{1}{\eps}(1+|\xi|^{-2})\||\xi|^2y^\alpha u\|_{L^p_m}\right)
\\[1.5ex]
&\leq  C\left(\|u-y^\alpha L_{2ta\cdot\xi}u+|\xi|^2y^\alpha u\|_{L^p_m}+\eps\|u\|_{W^{2,p}_{\mathcal N}(\alpha,m)}+\frac{1}{\eps}(1+|\xi|^{-2})\|u-y^\alpha L_{2ta\cdot\xi}u+|\xi|^2y^\alpha u\|_{L^p_m}\right).
\end{align*}
Note that for the last inequality we used the fact that the estimate in Lemma \ref{bound}  is uniform in $\xi$ and $a$.
By choosing $\eps=\frac{1}{2C}$ we deduce for some $C$ depending on $\xi$ but independent of $t$
\begin{align*}
&\|u\|_{W^{2,p}_{\mathcal N}(\alpha,m)}+\|y^\alpha u\|_{L^p_m}\leq C\|u-y^\alpha L_{2ta\cdot\xi}u+|\xi|^2y^\alpha u\|_{L^p_m}.
\end{align*}

Since, for $t=0$, the operator  $I-y^\alpha B+|\xi|^2y^\alpha $ is invertible in $W^{2,p}_{\mathcal N}(\alpha,m)\cap  L^p_{m+\alpha p}$ by Theorem \ref{1d}(ii),  the same holds for  $I-y^\alpha L_{2a\cdot\xi}+|\xi|^2y^\alpha$, by the method of continuity.
 
Let $(L_{m,p},D_{m,p})$ be the generator of $(e^{t(y^\alpha L_{2a\cdot\xi}-|\xi|^2y^\alpha)})$ in $L^p_m$ and consider the set $$\mathcal {D}= \left\{u \in C_c^\infty ([0,\infty)): u \ {\rm constant \ in \ a \ neighborhood\  of }\  0  \right \}$$ which is dense in  $W^{2,p}_{\mathcal N}(\alpha,m)\cap  L^p_{m+\alpha p}$, by Theorem \ref{core1-d}. 
By using the definition of $y^\alpha L_{2a\cdot\xi}-|\xi|^2y^\alpha$ through the form  $\mathfrak{a}$ as in the beginning of Section 3, it is easy to see that $\mathcal D \subset D_{c-\alpha,2}$ and that $y^\alpha L_{2a\cdot\xi}-|\xi|^2y^\alpha=L_{c-\alpha,2}$ on $\mathcal D$. Since $\mathcal D$ is dense in $W^{2,2}_{ \mathcal N}(\alpha,  c-\alpha)\cap L^2_{c+\alpha}$,  $y^\alpha L_{2a\cdot\xi}-|\xi|^2y^\alpha$ is closed on $W^{2,2}_{\mathcal N}(\alpha, c-\alpha)\cap L^2_{c+\alpha}$ and $L_{c-\alpha,2}$ is closed on $D_{c-\alpha,2}$, it follows that $W^{2,2}_{\mathcal N}(\alpha, c-\alpha)\cap L^2_{c+\alpha} \subset D_{c-\alpha,2}$ and then  $W^{2,2}_{\mathcal N}(\alpha,c-\alpha)\cap L^2_{c+\alpha} = D_{c-\alpha,2}$,  $y^\alpha L_{2a\cdot\xi}-|\xi|^2y^\alpha=L_{c-\alpha,2}$, since both operators are invertible on their own domains and one is an extension of the other. This completes the proof in the special case $p=2, m=c-\alpha$.

Take now $u \in \mathcal D$ and let $f=\lambda u-(y^\alpha L_{2a\cdot\xi}-|\xi|^2y^\alpha) u \in L^p_m \cap L^2_{c-\alpha}$ for large $\lambda$. Let $v \in D_{m,p}$ solve $\lambda v-L_{m,p}v=f$. Since the semigroups are consistent, $v$ coincides with the $L^2_{c-\alpha}$ solution which, by the previous step, is $u$. This gives $\mathcal D \subset D_{m,p}$ and that $y^\alpha L_{2a\cdot\xi}-|\xi|^2y^\alpha=L_{m,p}$ on $\mathcal D$ and, as before, one concludes the proof for $p<\infty$.
\qed

\smallskip
We remark that Proposition \ref{Gen DN} assures that $y^\alpha(L_{2a\cdot\xi}-|\xi|^2)$  generates a semigroup on $L^p_m$  under the milder assumption  $0<\frac{m+1}p<c+1-\alpha$. However, the hypothesis $(m+1)/p+\alpha>0$ must be added when $\alpha<0$ to have   $\mathcal D \subset L^p_{m+\alpha p}$.

\smallskip
%
As consequence we deduce the domain of the operator $y^\alpha L_{2a\cdot\xi}$ in the special case $\alpha=0$.

\begin{cor} \label{domain-alpha=0}
	Let    $1<p < \infty$, $0<\frac{m+1}{p}<c+1$ and $b\in\R$. Then the domain of the operator $L_{b}$ is $W^{2,p}_{\mathcal N}(0,m)$. The set
	$\mathcal {D}$ defined in \eqref{defD} 
	is a core.
\end{cor}
{\sc Proof.} By the arbitrariness of $\xi$, $a$ we can write  $b=2a\cdot\xi$. The required claim then follows from Theorem \ref{domainLb} since  the domains of   $L_{2a\cdot\xi}$ coincides with the one of  $L_{2a\cdot\xi}-|\xi|^2$  which for $\alpha=0$ is is $W^{2,p}_{\mathcal N}(0,m)\cap  L^p_{m}=W^{2,p}_{\mathcal N}(0,m)$.\qed

\smallskip 
 Using Corollary \ref{Cor Rb} (i) with $m$ replaced by $m-\alpha p$, we can characterize the domain of $L_{2a\cdot\xi}-\frac{\lambda}{y^\alpha}$. In what follows we write $D_{m,p}(A)$ to denote the domain of an operator $A$ on $L^p_m$.
\begin{cor} \label{domain1su}
Let  $1<p < \infty$, $\alpha^+<\frac{m+1}{p} <c+1$,  $\xi\in \R^N\setminus\{0\}$ and $\lambda\in C^+$. Then  the generator of $(e^{z (L_{2a\cdot\xi}-\frac{\lambda}{y^\alpha})})$ is the operator $L_{2a\cdot\xi}-\frac{\lambda}{y^\alpha}$ with  domain  $W^{2,p}_{\mathcal N}(0,m)\cap L^p_{m-\alpha p}$. In particular  the set
$\mathcal {D}$ defined in \eqref{defD} 
is a core.
\end{cor}
{\sc Proof.} We use (i) of Corollary \ref{Cor Rb} and  Theorem \ref{domainLb} with  $m$ replaced by $\tilde m :=m-\alpha p$ (note that the condition $\alpha^+<\frac{m+1}{p} <c+1$ and $\alpha^-<\frac{\tilde m+1}{p} <c+1-\alpha$ are equivalent) to obtain
\begin{align*}
	D_{m,p}\left(L_{2a\cdot\xi}-\frac{\lambda}{y^\alpha}\right)&=\left(|\xi|^2- L_{2a\cdot\xi}+\frac{\lambda}{y^\alpha}\right)^{-1}\left(L^p_m\right)=\left(\lambda-y^\alpha L_{2a\cdot\xi}+|\xi|^2 y^\alpha\right)^{-1}(L^p_{m-\alpha p})\\[1.5ex]
	&=D_{m-\alpha p,p}\left(y^\alpha L_{2a\cdot\xi}-|\xi|^2 y^\alpha\right)=W^{2,p}_{\mathcal N}(\alpha,m-\alpha p)\cap  L^p_{m}.
\end{align*}
Lemma \ref{Sobolev Formula magica} then implies
\begin{align*}
	D_{m,p}\left(L_{2a\cdot\xi}-\frac{\lambda}{y^\alpha}\right)&=W^{2,p}_{\mathcal N}(0,m)\cap L^p_{m-\alpha p}.
\end{align*}

\qed

\begin{os}
Let   $1<p < \infty$ such that  $\alpha^{-}<\frac{m+1}{p}<c+1-\alpha$,  $\xi\in \R^N\setminus\{0\}$ and $\lambda\in\C^+$. Theorem \ref{domainLb}, Corollary \ref{domain1su} and   Lemma \ref{Sobolev Formula magica} show that the operators
\begin{align*}
y^\alpha L_{2a\cdot\xi}-|\xi|^2 y^\alpha\quad \text{in\quad $L^p_m$},\qquad\quad  L_{2a\cdot\xi}-\frac{\lambda}{y^\alpha}\quad \text{in \quad  $L^p_{m+\alpha p}$}
\end{align*}
 endowed with the common domain
\begin{align*}
	W^{2,p}_{\mathcal N}(\alpha,m)\cap L^p_{m+\alpha p}=W^{2,p}_{\mathcal N}(0,m+\alpha p)\cap L^p_m
\end{align*}
generate a semigroup on $L^p_m$  and $L^p_{m+\alpha p}$, respectively.
 Their resolvents satisfy
\begin{align*}
	\left(\lambda-y^\alpha L_{2a\cdot\xi}+|\xi|^2 y^\alpha\right)^{-1}f=\left(|\xi|^2- L_{2a\cdot\xi}+\frac{\lambda}{y^\alpha}\right)^{-1}\left(\frac{f}{y^\alpha}\right),\qquad f\in L^p_m
\end{align*}
proving the equivalence between the two elliptic equations
\begin{align*}
	\lambda u-y^\alpha L_{2a\cdot\xi}u+|\xi|^2 y^\alpha u=f,\qquad \quad 
	\frac{\lambda}{y^\alpha} u- L_{2a\cdot\xi}u+|\xi|^2  u=\frac{f}{y^\alpha}.	
\end{align*}
\end{os}

\section{Multipliers} \label{mult}
In this section we investigate the boundedness of some multipliers related to  the   degenerate operator
\begin{align} \label{Elle}
	\mathcal L =y^\alpha \left(\Delta_{x} +2\sum_{i=1}^Na_{i}D_{iy}+D_{yy}+\frac{c}{y}D_y\right) ,\qquad a \in\R^N,\ |a|<1,\qquad \alpha<2
\end{align}
Assuming that
$$y^\alpha\left(\Delta_x u+2a\cdot \nabla_xD_yu+ B_yu\right)=f$$
and taking the Fourier transform (denoted by $\mathcal F$ or $\hat \cdot$)  with respect to $x$ (with covariable $\xi$) we obtain 
$$ -y^\alpha|\xi|^2\,\hat u(\xi,y)+y^\alpha i2a\cdot\xi  D_y\hat u(\xi,y)+  y^\alpha B _y \hat u(\xi,y)=\hat f(\xi,y).$$ 
We  
  consider the operator $L_{2a\cdot \xi}=B+2ia\cdot \xi D_y$ of Section \ref{section L_b}. The latter computation  shows that formally 
\begin{align*}
	(\lambda-\mathcal L)^{-1}={\cal F}^{-1}\left (\lambda- y^\alpha L_{2a\cdot\xi}+  y^\alpha|\xi|^2 \right)^{-1} {\cal F}.
\end{align*} 
In order to prove that $\mathcal L$ generates an analytic semigroup and to prove regularity for the associated parabolic problem, we   investigate  the boundedness of the  operator-valued multiplier
\begin{align*}
	\xi \in \R^N \to R_{\lambda}(\xi)&= \left (\lambda- y^\alpha L_{2a\cdot\xi}+   y^\alpha|\xi|^2  \right)^{-1}.
\end{align*}
To characterize the  domain of $\mathcal L$ we also consider the multipliers   $|\xi|^2y^\alpha R_\lambda$, $\xi y^\alpha D_yR_{\lambda}$ which are associated  with the operators $	y^\alpha\Delta_x(\lambda-\mathcal L)^{-1}$, $ y^\alpha D_{xy}	(\lambda-\mathcal L)^{-1}$, respectively. In the next results  we  prove that the  above  multipliers  satisfy the hypotheses of Theorem \ref{marcinkiewicz}.

\smallskip 
We also need  the operator-valued multiplier $\tilde R_{\lambda}(\xi)$ defined by 
\begin{align*}
\xi \in \R^N \to \tilde R_{\lambda}(\xi)&= \left (|\xi|^2 -L_{2a\cdot\xi}+\frac{\lambda}{y^\alpha}  \right)^{-1}\in \mathcal B\left(L^p_{m+\alpha p}\right).
\end{align*} 
Note that the role of $\xi$ and $\lambda$ is interchanged between $R_\lambda (\xi)$ and $\tilde R_\lambda (\xi)$: $\lambda $ is a spectral parameter in the first and plays a role of a complex potential in the second, where instead is $\xi$  the spectral parameter.  Nevertheless, we keep the same notation.

By  Corollary \ref{Cor Rb} we have
\begin{align}\label{Formule Magiche}
	y^\alpha R_\lambda(\xi)f&=y^\alpha\tilde R_{\lambda}(\xi)\left(\frac{f}{y^\alpha}\right),\qquad y^\alpha D_y R_\lambda(\xi)f=y^\alpha D_y \tilde R_{\lambda}(\xi)\left(\frac{f}{y^\alpha}\right),\qquad f\in L^p_m.
\end{align}

\begin{prop}\label{Mikhlin 0order}
Let $1<p<\infty$, $\alpha<2$ be such that $\alpha^-<\frac{m+1}p<c+1-\alpha$. 
Then  the families  
\begin{align*}
	\left\{\lambda  R_{\lambda}(\xi),\; |\xi|^2y^\alpha R_{\lambda}(\xi),\; \xi y^\alpha D_y  R_\lambda(\xi)\in \mathcal{B}(L^p_{m})\,:\, \lambda\in \C_+,\; \xi\in\R^N\setminus\{0\}\right\}
\end{align*}
are $\mathcal R$-bounded.
\end{prop}
{\sc Proof.} The $\mathcal R$-boundedness of $\lambda R_\lambda(\xi)$ follows by  Corollary \ref{Rbound lR_l}. 

 The $\mathcal R$-boundedness of $|\xi|^2 y^\alpha  R_\lambda(\xi)$ in $\mathcal{B}(L^p_{m})$ follows by using formula \eqref{Formule Magiche}. We write
 \begin{align*}
 	|\xi|^2	y^\alpha R_\lambda(\xi)&=y^\alpha\left(|\xi|^2\tilde R_{\lambda}(\xi)\right)\left(\frac{\cdot}{y^\alpha}\right)
 \end{align*} 
and use  the $\mathcal R$-boundedness of $|\xi|^2 \tilde R_\lambda(\xi)$  in   $\mathcal{B}(L^p_{m+\alpha p})$ proved  in  Proposition \ref{kernel Vreal} with $V(y)=\lambda y^{-\alpha}$.

Let us finally prove the $\mathcal R$-boundedness of $\xi y^\alpha  D_y R_\lambda(\xi)$ in $\mathcal{B}(L^p_{m})$. By   formula \eqref{Formule Magiche} again we have
\begin{align*}
	\xi y^\alpha D_yR_{\lambda}(\xi)= y^\alpha \left(\xi D_y\tilde R_{\lambda}(\xi)\right)\left(\frac{\cdot }{y^\alpha}\right).
\end{align*}
Let us write
\begin{align*}
	\tilde R_{\lambda}(\xi)&= \left (|\xi|^2 -L_{2a\cdot\xi}+\frac{\lambda}{y^\alpha}  \right)^{-1}\\[1ex]
	&=
	\left (|\xi|^2 -L_{2a\cdot\xi}\right)^{-1}\left (|\xi|^2 -L_{2a\cdot\xi}\right)\left (|\xi|^2 -L_{2a\cdot\xi}+\frac{\lambda}{y^\alpha}  \right)^{-1}\\[1ex]
	&=\left (|\xi|^2 -L_{2a\cdot\xi}\right)^{-1}\left(Id-\frac{\lambda}{y^\alpha}\tilde R_{\lambda}(\xi)\right)\\[1ex]
	&=\left (|\xi|^2 -L_{2a\cdot\xi}\right)^{-1}\left(\frac{\cdot}{y^\alpha }\right)\left(y^\alpha-\lambda \tilde R_{\lambda}(\xi)\right)
\end{align*} 
(note that   $\left (|\xi|^2 -L_{2a\cdot\xi}\right)\left (|\xi|^2 -L_{2a\cdot\xi}+\frac{\lambda}{y^\alpha}  \right)^{-1}$ is well defined since  $D\left(L_{2a\cdot\xi}-\frac{\lambda}{y^\alpha}\right)\subseteq D\left(L_{2a\cdot\xi}\right)$, by  Corollary \ref{domain1su}).
The previous relations and \eqref{Formule Magiche} give
\begin{align*}
	\xi y^\alpha D_y R_{\lambda}(\xi)&=\xi y^\alpha  D_y\left (|\xi|^2 -L_{2a\cdot\xi}\right)^{-1}\left(\frac{\cdot }{y^\alpha}\right) \left(y^\alpha-\lambda \tilde R_{\lambda}(\xi)\right)\left(\frac{\cdot }{y^\alpha}\right)\\[1.5ex]
	&=\xi y^\alpha  D_y\left (|\xi|^2 -L_{2a\cdot\xi}\right)^{-1}\left(\frac{\cdot }{y^\alpha}\right) \left(I-\lambda \tilde R_{\lambda}(\xi)\left(\frac{\cdot }{y^\alpha}\right)\right)\\[1.5ex]
	&= \xi y^\alpha D_y\left (|\xi|^2 -L_{2a\cdot\xi}\right)^{-1}\left(\frac{\cdot }{y^\alpha}\right) \left(I- \lambda R_{\lambda}(\xi)\right).
\end{align*}
The $\mathcal R$-boundedness of the family $\xi y^\alpha  D_y R_\lambda(\xi)$ in $\mathcal{B}(L^p_{m})$ then follows by composing  the   $\mathcal R$-boundedness of $\xi D_y\left (|\xi|^2 -L_{2a\cdot\xi}\right)^{-1}$ in $\mathcal{B}(L^p_{m+\alpha p})$ proved in \cite[Corollary 6.4]{MNS-Singular-Half-Space} and the $\mathcal R$-boundedness of the family $\lambda  R_\lambda(\xi)$ in $\mathcal{B}(L^p_{m})$ proved in the first step.
\qed

To apply the Mikhlin multiplier theorem, we need a formula for the derivatives  of the above functions with respect to $\xi$. In the following lemma ${\mathcal S}_n$ denotes the set of permutations of $n$ elements.

\begin{lem}   \label{Lema Marck1}
Let $1<p<\infty$ , $\alpha<2$ be such that $\alpha^-<\frac{m+1}p<c+1-\alpha$, and let   us consider, for any fixed $\lambda\in\C_+$, the  map  
\begin{align*}
	\xi \in \R^N \to R_{\lambda}(\xi)&= (\lambda-y^\alpha L_{2a\cdot\xi}+y^\alpha |\xi|^2)^{-1}\in B(L^p_m).
\end{align*} 
Then $R_\lambda,\;y^\alpha R_\lambda,\; y^\alpha D_yR_\lambda\in C^\infty\left(\R^N\setminus \{0\}; B(L^p_m)\right)$ 
and for any family of different  indexes $j_1,j_2,\dots,j_n\in \{1,\dots,N\}$ one has   \begin{align}\label{Explicit Derivative Rlamba}
	\nonumber 	D_{\xi_{j_1}}\cdots D_{\xi_{j_n}} R_\lambda(\xi)&=\sum_{\sigma\in {\mathcal S}_n} R_\lambda(\xi) \prod_{k=1}^n\Big( 2ia_{j_{\sigma(k)}}y^\alpha D_yR_\lambda(\xi)-2\xi_{j_{\sigma(k)}}y^\alpha R_\lambda(\xi)\Big)\\[1ex]
		D_{\xi_{j_1}}\cdots D_{\xi_{j_n}} y^\alpha R_\lambda(\xi)&=\sum_{\sigma\in {\mathcal S}_n} y^\alpha R_\lambda(\xi) \prod_{k=1}^n\Big( 2ia_{j_{\sigma(k)}}y^\alpha D_yR_\lambda(\xi)-2\xi_{j_{\sigma(k)}}y^\alpha R_\lambda(\xi)\Big)\\[1ex]\nonumber 
	D_{\xi_{j_1}}\cdots D_{\xi_{j_n}} y^\alpha D_yR_\lambda(\xi)&=\sum_{\sigma\in {\mathcal S}_n} y^\alpha D_y R_\lambda(\xi) \prod_{k=1}^n\Big( 2ia_{j_{\sigma(k)}}y^\alpha D_yR_\lambda(\xi)-2\xi_{j_{\sigma(k)}}y^\alpha R_\lambda(\xi)\Big).
\end{align}
\end{lem}
{\sc{Proof.}}
Let us fix  $\lambda\in \C_+$. Let us prove the first equality in  \eqref{Explicit Derivative Rlamba} for $n=1$  that is, for $j=1,\dots, n$   
\begin{align}\label{Lema Marck1 eq2}
\frac{\partial}{\partial \xi_j} (R_{\lambda}(\xi))= R_{\lambda}(\xi) \Big(2ia_jy^\alpha D_y R_{\lambda}(\xi)-2\xi_jy^\alpha  R_{\lambda}(\xi) \Big),\qquad \xi\in\R^n\setminus\{0\}.
\end{align}
Indeed let us  write for $|h|\leq 1$
\begin{align}\label{Lema Marck1 eq1}
\nonumber  R_{\lambda}&(\xi+he_j)-R_{\lambda}(\xi)= \left(\lambda+y^\alpha |\xi+he_j|^2-y^\alpha L_{2a\cdot (\xi+he_j)}\right)^{-1}-\left(\lambda+y^\alpha |\xi|^2-y^\alpha L_{2a\cdot \xi}\right)^{-1}\\[2ex]
\nonumber&= R_{\lambda}(\xi)\,\Big[\left(\lambda+y^\alpha |\xi|^2-y^\alpha L_{2a\cdot \xi}\right)\left(\lambda+y^\alpha |\xi+he_j|^2-y^\alpha L_{2a\cdot (\xi+he_j)}\right)^{-1}-I\Big]\\[2ex]
\nonumber&= R_{\lambda}(\xi)\,y^\alpha \left( L_{2a\cdot (\xi+he_j)}-  L_{2a\cdot \xi}+|\xi|^2-|\xi+he_j|^2\right)\left(\lambda+y^\alpha |\xi+he_j|^2-y^\alpha L_{2a\cdot (\xi+he_j)}\right)^{-1}\\[2ex]\nonumber 
&=R_{\lambda}(\xi)\,y^\alpha \left(2ia_jhD_y-2\xi_jh-h^2\right)R_{\lambda}(\xi+he_j)\\[2ex]
&=2ia_jh\, R_{\lambda}(\xi)\,y^\alpha D_yR_{\lambda}(\xi+he_j)- (2\xi_jh+h^2)\,R_{\lambda}(\xi)\,y^\alpha R_{\lambda}(\xi+he_j).
\end{align}
Using the $\mathcal{R}$-boundedness of  $\lambda  R_{\lambda}(\xi)$ $|\xi|^2y^\alpha R_{\lambda}(\xi)$, $\xi y^\alpha D_yR_{\lambda}(\xi)$ of Proposition \ref{Mikhlin 0order} (which implies uniform boundedness), the last equation  implies in particular that
\begin{align}\label{Lema Marck1 eq5}\nonumber 
	&R_{\lambda}(\xi+he_j)\to R_{\lambda}(\xi),\qquad\qquad y^\alpha R_{\lambda}(\xi+he_j)\to y^\alpha  R_{\lambda}(\xi),\\[2ex]
	&y^\alpha D_y R_{\lambda}(\xi+he_j)\to y^\alpha D_y R_{\lambda}(\xi)\qquad \qquad \text{in the norm of $B\left(L^p_m\right)$ \quad as $h\to 0$.} 
\end{align}
 For example  from \eqref{Lema Marck1 eq1} one has  for some positive constant $C$
\begin{align*} 
 &\left\|R_{\lambda}(\xi+he_j)-R_{\lambda}(\xi)\right\|_{\mathcal{B}(L^p_m)}\leq C \frac{\left\|R_\lambda(\xi)\right\|_{\mathcal{R}\left(\mathcal{B}(L^p_m)\right)}}{|\lambda|}\\[1ex]
 &\hspace{10ex}\times \left(|h|\frac{\left\|y^\alpha D_y R_\lambda(\xi)\right\|_{\mathcal{R}\left(\mathcal{B}(L^p_m)\right)}}{|\xi+he_j|}+(|\xi|h+h^2)\frac{\left\|y^\alpha  R_\lambda(\xi)\right\|_{\mathcal{R}\left(\mathcal{B}(L^p_m)\right)}}{|\xi+he_j|^2}\right)\to 0\qquad \text{as}\quad h\to 0.
\end{align*}
The other limits in \eqref{Lema Marck1 eq5} follow similarly after applying to  both sides of \eqref{Lema Marck1 eq1}   $y^\alpha $ and $y^\alpha D_y$, respectively.

To end the proof we apply equality \eqref{Lema Marck1 eq1} again  to get
\begin{align}\label{Rapporto incrementale}
	\frac{R_{\lambda}(\xi+he_j)-R_{\lambda}(\xi)}h&=R_{\lambda}(\xi)y^\alpha  \left(2ia_jD_y-2\xi_j\right) R_{\lambda}(\xi)- h\,R_{\lambda}(\xi)\,y^\alpha R_{\lambda}(\xi+he_j)
\end{align}
which tends to  $R_{\lambda}(\xi)y^\alpha  \left(2ia_jD_y-2\xi_j\right) R_{\lambda}(\xi)$ in the norm of $B\left(L^p_m\right)$  as $h\to 0$ since , by \eqref{Lema Marck1 eq5}, the last term  tends to $0$.  This proves \eqref{Lema Marck1 eq2}.

The proof of  the other equalities in  \eqref{Explicit Derivative Rlamba} for $n=1$  that is, for $j=1,\dots, n$ , 
\begin{align*}
	\frac{\partial}{\partial \xi_j} (y^\alpha R_{\lambda}(\xi))= y^\alpha R_{\lambda}(\xi) \left(2ia_j y^\alpha D_y R_{\lambda}(\xi)-2\xi_j y^\alpha R_{\lambda}(\xi)\right)  ,\qquad \xi\in\R^n\setminus\{0\},\\[1ex]
	\frac{\partial}{\partial \xi_j} (y^\alpha D_yR_{\lambda}(\xi))= y^\alpha D_yR_{\lambda}(\xi) \left(2ia_j y^\alpha D_y R_{\lambda}(\xi)-2\xi_j y^\alpha R_{\lambda}(\xi)\right)  ,\qquad \xi\in\R^n\setminus\{0\}
\end{align*}
follow similarly by applying, respectively,  the operators $y^\alpha Id$, $y^\alpha D_y$ to both sides of   \eqref{Rapporto incrementale}  and taking the limit for $h\to 0$. For example  for  the derivative $\frac{\partial}{\partial \xi_j} (y^\alpha D_yR_{\lambda}(\xi))$  we write,  as in \eqref{Lema Marck1 eq1},
\begin{align*}
	\frac{y^\alpha D_y R_{\lambda}(\xi+he_j)-y^\alpha D_yR_{\lambda}(\xi)}h=&y^\alpha D_y R_{\lambda}(\xi)y^\alpha  \left(2ia_jD_y-2\xi_j\right) R_{\lambda}(\xi)\\[1ex]
	&- h\,y^\alpha D_y R_{\lambda}(\xi)\,y^\alpha R_{\lambda}(\xi+he_j)
\end{align*}
which by \eqref{Lema Marck1 eq5} again  tends to  $y^\alpha D_y R_{\lambda}(\xi)y^\alpha  \left(2ia_jD_y-2\xi_j\right) R_{\lambda}(\xi)$ in the norm of $B\left(L^p_m\right)$  as $h\to 0$.

Finally,  \eqref{Explicit Derivative Rlamba} for $n>1$  follows by induction. For example if $n=2$ and $l\neq j$ one has 
\begin{align*}
	\frac{\partial^2}{\partial \xi_l\partial \xi_j} (R_{\lambda}(\xi))&=\frac{\partial}{\partial \xi_l}\Big[ R_{\lambda}(\xi) \Big(2ia_jy^\alpha D_y R_\lambda(\xi)-2\xi_jy^\alpha R_{\lambda}(\xi)\Big) \Big]\\[2ex]
	&=\frac{\partial}{\partial \xi_l}(R_{\lambda}(\xi))\Big(2ia_jy^\alpha D_y R_\lambda(\xi)-2\xi_jy^\alpha R_{\lambda}(\xi)\Big)\\[1ex]
	&\hspace{2ex}+ R_{\lambda}(\xi)\left(2ia_j\frac{\partial}{\partial \xi_l}\left(y^\alpha D_y R_\lambda(\xi)\right)-2\xi_j\frac{\partial}{\partial \xi_l}\left(y^\alpha R_{\lambda}(\xi)\right)\right)\\[2ex]
	&=R_{\lambda}(\xi) \Big(2ia_ly^\alpha D_y R_\lambda(\xi)-2\xi_l y^\alpha R_\lambda(\xi)\Big)\Big(2ia_jy^\alpha D_y R_\lambda(\xi)-2\xi_jy^\alpha R_{\lambda}(\xi)\Big)\\[2ex]
	&\hspace{2ex}+ R_{\lambda}(\xi) \Big(2ia_jy^\alpha D_y R_\lambda(\xi)-2\xi_jy^\alpha R_{\lambda}(\xi)\Big) \Big(2ia_ly^\alpha D_y R_\lambda(\xi)-2\xi_l y^\alpha R_\lambda(\xi)\Big).
\end{align*}
\qed

Now we can finally prove that the multiplier  $\lambda R_\lambda(\xi)$ associated with the operators  $\lambda(\lambda-\mathcal L)^{-1}$ satisfies the hypothesis of  Theorem \ref{marcinkiewicz}. This is crucial    for proving   that $	\mathcal L $ generates an analytic semigroup  in $L^p_m$.

\begin{teo}\label{Teo Mult Resolv}
Let $1<p<\infty$ , $\alpha<2$ be such that $\alpha^-<\frac{m+1}p<c+1-\alpha$.
Then the family 
$$\left \{\xi^{\beta}D^\beta_\xi(\lambda  R_{\lambda}(\xi)): \xi\in \R^{N}\setminus\{0\}, \ \beta  \in \{0,1\}^N ,\lambda \in \C_+\right \}$$
is $\mathcal{R}$-bounded in $L^p_m$.
\end{teo}
{\sc{Proof.}} Let $\beta=(\beta_1,\dots,\beta_N)\in\{0,1\}^N$, $|\beta|=n$.  Let us suppose, without any loss of generality, $\beta_i=1$, for $i\leq n$ and  $\beta_i=0$, for $i>n$. Then using  \eqref{Explicit Derivative Rlamba} we get
\begin{align*}	
	\xi^\beta D_\xi^\beta \lambda R_\lambda(\xi)&=\xi_1\cdots \xi_n \, \left(D_{\xi_{1}}\cdots D_{\xi_{n}}\right) \lambda R_\lambda(\xi)\\[1ex]
	&=\sum_{\sigma\in {\mathcal S}_n} \lambda R_\lambda(\xi) \prod_{j=1}^n\Big( 2ia_{\sigma(j)}\xi_{\sigma(j)}y^\alpha D_yR_\lambda(\xi)-2\xi_{\sigma(j)}^2y^\alpha R_\lambda(\xi)\Big).
\end{align*} 

The  $\mathcal R$-boundedness of $\xi^{\beta}D^\beta_\xi(\lambda  R_{\lambda})(\xi)$ then follows by composition and domination from the $\mathcal R$-boundedness of $\lambda R_{\lambda}(\xi)$, $|\xi|^2y^\alpha R_{\lambda}(\xi)$, $\xi y^\alpha R_{\lambda}(\xi)$ using  Proposition \ref{Mikhlin 0order} and  Corollary \ref{domination}.
\qed
\medskip 

%

The next two theorems  show  that the multipliers  $|\xi|^2 y^\alpha R_\lambda$, $\xi y^\alpha D_yR_\lambda$, associated respectively  with the operators  $y^\alpha \Delta_x(\lambda-\mathcal L)^{-1}$, $y^\alpha D_{xy}(\lambda-\mathcal L)^{-1}$, satisfy the hypotheses of   Theorem \ref{marcinkiewicz}. This is essential     for characterizing the domain of  $\mathcal L $.
\begin{teo}\label{Mult xi^2 Resolv}
Let $1<p<\infty$ , $\alpha<2$ be such that $\alpha^-<\frac{m+1}p<c+1-\alpha$. 
Then   the family 
$$\left \{\xi^{\beta}D^\beta_\xi(|\xi|^2 y^\alpha R_{\lambda}(\xi)): \xi\in \R^{N}\setminus\{0\}, \ \beta  \in \{0,1\}^N ,\lambda \in \C_+\right \}$$
is $\mathcal{R}$-bounded in $L^p_m$.
\end{teo}
{\sc{Proof.}}  Let us prove preliminarily that  the family 
$$\left \{|\xi|^2\xi^{\beta}D^\beta_\xi( y^\alpha R_{\lambda}(\xi)): \xi\in \R^{N}\setminus\{0\}, \ \beta  \in \{0,1\}^N ,\lambda \in \C_+\right \}$$
is $\mathcal{R}$-bounded in $L^p_m$.  Let $\beta=(\beta_1,\dots,\beta_N)\in\{0,1\}^N$, $|\beta|=n$. Let us suppose, without any loss of generality, $\beta_i=1$, for $i\leq n$ and  $\beta_i=0$, for $i>n$. Then using \eqref{Explicit Derivative Rlamba} one has 
\begin{align*}	
	|\xi|^2\xi^\beta D_\xi^\beta y^\alpha  R_\lambda(\xi)&=	|\xi|^2\xi_1\cdots \xi_n \, \left(D_{\xi_{1}}\cdots D_{\xi_{n}}\right) y^\alpha R_\lambda(\xi)\\[1ex]
	&=\sum_{\sigma\in {\mathcal S}_n} |\xi|^2y^\alpha  R_\lambda(\xi) \prod_{j=1}^n\Big( 2ia_{\sigma(j)}\xi_{\sigma(j)}y^\alpha D_yR_\lambda(\xi)-2\xi_{\sigma(j)}^2y^\alpha R_\lambda(\xi)\Big)
\end{align*} 
and the required  $\mathcal R$-boundedness of $	|\xi|^2\xi^\beta D_\xi^\beta y^\alpha  R_\lambda(\xi)$ then follows as at  the end  of Theorem \ref{Teo Mult Resolv}.

To prove the required claim let us observe that for any  $\beta  \in \{0,1\}^N$, $|\beta|=n$  there exist $\beta^j  \in \{0,1\}^N$ satisfying
\begin{align*}
	\beta^j_k=\beta_k,\quad k\neq j,\qquad  \beta^j_j=0,\qquad |\beta^j|=n-1
\end{align*}
and  such that
\begin{align*}
D^\beta_\xi(|\xi|^2  y^\alpha R_{\lambda}(\xi))=\sum_{j:\beta_j=1}2\xi_jD^{\beta^j}_\xi y^\alpha  R_{\lambda}(\xi)+|\xi|^2D^\beta_\xi y^\alpha R_{\lambda}(\xi). 
\end{align*}
Then 
\begin{align*}
\xi^\beta	D^\beta_\xi(|\xi|^2  y^\alpha R_{\lambda}(\xi))=\sum_{j:\beta_j=1}2 \xi^2_j\xi^{\beta^j}D^{\beta^j}_\xi y^\alpha  R_{\lambda}(\xi)+|\xi|^2\xi^\beta D^\beta_\xi y^\alpha R_{\lambda}(\xi).
\end{align*}
and the proof now follows by domination using the previous step.\qed

\begin{teo}\label{Mult D_y/y Resolv}
Let $1<p<\infty$ , $\alpha<2$ be such that $\alpha^-<\frac{m+1}p<c+1-\alpha$. 
Then the family 
$$\left \{\xi^{\beta}D^\beta_\xi\left( \xi y^\alpha D_yR_{\lambda}(\xi)\right): \xi\in \R^{N}\setminus\{0\}, \ \beta  \in \{0,1\}^N ,\lambda \in \C_+\right \}$$
is $\mathcal{R}$-bounded in $L^p_m$.
\end{teo}
{\sc{Proof.}}  As in the  proof  of the previous theorem we prove preliminarily that  the family 
$$\left \{\xi \,\xi^{\beta}D^\beta_\xi( y^\alpha D_y R_{\lambda}(\xi)): \xi\in \R^{N}\setminus\{0\}, \ \beta  \in \{0,1\}^N ,\lambda \in \C_+\right \}$$
is $\mathcal{R}$-bounded in $L^p_m$. Indeed let $\beta=(\beta_1,\dots,\beta_N)\in\{0,1\}^N$, $|\beta|=n$. Let us suppose, without any loss of generality, $\beta_i=1$, for $i\leq n$ and  $\beta_i=0$, for $i>n$. Then using \eqref{Explicit Derivative Rlamba} one has 
\begin{align*}	
	\xi\, \xi^\beta D_\xi^\beta y^\alpha  D_yR_\lambda(\xi)&=	\xi\, \xi_1\cdots \xi_n \, \left(D_{\xi_{1}}\cdots D_{\xi_{n}}\right) y^\alpha D_yR_\lambda(\xi)\\[1ex]
	&=\sum_{\sigma\in {\mathcal S}_n} \xi\,y^\alpha  D_y R_\lambda(\xi) \prod_{j=1}^n\Big( 2ia_{\sigma(j)}\xi_{\sigma(j)}y^\alpha D_yR_\lambda(\xi)-2\xi_{\sigma(j)}^2y^\alpha R_\lambda(\xi)\Big)
\end{align*} 
and the required  $\mathcal R$-boundedness of $	\xi\,\xi^\beta D_\xi^\beta y^\alpha  D_yR_\lambda(\xi)$ then follows as at  the end  of Theorem \ref{Teo Mult Resolv}.

To prove the required claim let us fix   $\beta  \in \{0,1\}^N$, $|\beta|=n$ and let us observe that for any $j=1,\dots N$ one has  
\begin{align*}
D^\beta_\xi(\xi_{j} y^\alpha D_yR_{\lambda})(\xi)=\xi_{j}D^{\beta}_\xi \left(y^\alpha D_yR_{\lambda}(\xi)\right)+\beta_j D^{\beta^j}_\xi \left(y^\alpha D_yR_{\lambda}(\xi)\right). 
\end{align*}
  where $\beta^j  \in \{0,1\}^N$ satisfies
\begin{align*}
	\beta^j_k=\beta_k,\quad k\neq j,\qquad  \beta^j_j=0,\qquad |\beta^j|=n-1.
\end{align*}
The proof then follows by the previous step.
\qed

\section{The operator  $\mathcal L=y^\alpha(\Delta_{x} +2a\cdot\nabla_xD_y+ B_y), \quad |a|<1, \alpha<2$} \label{Sect DOm max}

In this section we prove generation results, maximal regularity and   domain characterization for  the    operator $\mathcal L$ defined in (\ref{Elle})
 %
in $L^p_m$.  
More general operators will be treated in the next section, based on this model case. We start with the $L^2$ theory.

 As explained at the beginning of Section \ref{mult}, we have formally for $\lambda\in\C_+$
\begin{align*}
	\lambda(\lambda-\mathcal L)^{-1}&={\cal F}^{-1}\left (\lambda R_\lambda(\xi)\right) {\cal F},\qquad  R_{\lambda}(\xi)= \left (\lambda- y^\alpha L_{2a\cdot\xi}+   y^\alpha|\xi|^2  \right)^{-1}
\end{align*} 
and consequently 
\begin{align*}
	y^\alpha \Delta_x (\lambda-\mathcal L)^{-1}&=-{\cal F}^{-1} \left(|\xi|^2y^\alpha R_\lambda (\xi)\right) {\cal F},\qquad  i\nabla_x D_y (\lambda-\mathcal L)^{-1}={\cal F}^{-1} \left(\xi y^\alpha D_yR_\lambda (\xi)\right) {\cal F}.
\end{align*} 
All  properties of  $\mathcal L$ follow from   the boundedness of the above multipliers,    through Theorem \ref{marcinkiewicz}.

\subsection{The operator $\mathcal L$ in $L^2_{c-\alpha}$} 
We  assume that $c+1-\alpha>0$ so that the measure $y^{c-\alpha}\, dx\, dy$ is locally finite near $y=0$ and use the Sobolev space $H^{1}_{\alpha, c}:=\{u \in L^2_{c-\alpha} : y^\frac{\alpha}{2}\nabla u \in L^2_{c-\alpha}\}$ equipped with the inner product
\begin{align*}
	\left\langle u, v\right\rangle_{H^1_{\alpha, c}}:= \left\langle u, v\right\rangle_{L^2_{c-\alpha}}+\left\langle y^\frac{\alpha}{2}\nabla u, y^\frac{\alpha}{2}\nabla v\right\rangle_{L^2_{c-\alpha}}.
\end{align*}
We consider the  form  in $L^2_{c-\alpha}$ 
\begin{align*}
	\mathfrak{a}(u,v)
	&:=
	\int_{\R^{N+1}_+} \langle \nabla u, \nabla \overline{v}\rangle\,y^{c} dx\,dy+2\int_{\R^{N+1}_+} D_yu\, a\cdot\nabla_x\overline{v}\,y^{c} dx\,dy , \quad
	D(\mathfrak{a})=H^1_{\alpha, c}
\end{align*}

and its adjoint $\mathfrak{a}^*(u,v)=\overline{\mathfrak{a}(v,u)}$  
\begin{align*}
	\mathfrak{a^*}(u,v)=\overline{\mathfrak{a}(v,u)}
	&:=
	\int_{\R^{N+1}_+} \langle \nabla u, \nabla \overline{v}\rangle\,y^{c} dx\,dy+2\int_{\R^{N+1}_+} a\cdot \nabla_x u\,D_y\overline{v}\,y^{c} dx\,dy.
\end{align*}

\begin{prop} \label{prop-form}
The forms $\mathfrak{a}$, $\mathfrak{a^*}$ are continuous, accretive and   sectorial.
\end{prop}
{\sc Proof.} We consider only the form $\mathfrak{a}$, the adjoint form can be handled similarly. If $u\in H^1_{\alpha, c}$
\begin{align*}
\Rp \mathfrak{a}(u,u)\geq \|\nabla_x u\|^2_{L^2_c}+\|D_y u\|^2_{L^2_c}-2|a| \|\nabla_x u\|_{L^2_c}\|D_y u\|_{L^2_c}\geq (1-|a|)(\|\nabla_x u\|^2_{L^2_c}+\|D_y u\|^2_{L^2_c}).
\end{align*}
By the ellipticity assumption $|a|<1$, the accretivity follows.
Moreover 
\begin{align*}
|\Ip \mathfrak{a}(u,u)|\leq  2|a|\|\nabla_x u\|_{L^2_c}\|D_y u\|_{L^2_c}\leq |a|(\|\nabla_x u\|^2_{L^2_c}+\|D_y u\|^2_{L^2_c})\leq \frac{|a|}{(1-|a|)}\Rp \mathfrak{a}(u,u).
\end{align*}
This proves the sectoriality and then the continuity of the form. \qed

We define the operators $\mathcal L$ and $\mathcal L^*$ associated respectively to the forms $\mathfrak{a}$ and $\mathfrak{a}^*$  by
\begin{align} \label{BesselN}
	\nonumber D( \mathcal L)&=\{u \in H^1_{\alpha, c}: \exists  f \in L^2_{c-\alpha} \ {\rm such\ that}\  \mathfrak{a}(u,v)=\int_{\R^{N+1}_+} f \overline{v}y^{c}\, dz\ {\rm for\ every}\ v\in H^1_{\alpha, c}\},\\  \mathcal Lu&=-f;
\end{align}
\begin{align} \label{adjoint}
	\nonumber D( \mathcal L^*)&=\{u \in H^1_{\alpha, c}: \exists  f \in L^2_{c-\alpha} \ {\rm such\ that}\  \mathfrak{a}^*(u,v)=\int_{\R^{N+1}_+} f \overline{v}y^{c}\, dz\ {\rm for\ every}\ v\in H^1_{\alpha, c}\},\\  \mathcal L^*u&=-f.
\end{align}
If $u,v$ are smooth function with compact support in the closure of $\R_+^{N+1}$ (so that they do not need to vanish on the boundary), it is easy to see integrating by parts that $$-\mathfrak a (u,v)= \langle y^\alpha(\Delta_x u+2a\cdot \nabla_xD_yu+ B_yu), \overline v\rangle_{L^2_{c-\alpha}}  $$
if $\ds \lim_{y \to 0} y^c D_yu(x,y)=0$. This means that 
$\mathcal L$ is the operator $y^\alpha(\Delta_x +2a\cdot \nabla_xD_y+ B_y)$ with Neumann boundary conditions at $y=0$. On the other hand 
$$-\mathfrak a^* (u,v)= \left\langle y^\alpha\left(\Delta_x u+2a\cdot \nabla_xD_yu+2(c-\alpha)\frac {a\cdot \nabla_x u}{y}+ B_yu\right), \overline v\right\rangle_{L^2_{c-\alpha}}  $$
if $\ds \lim_{y \to 0} y^{c} \left ( D_yu(x,y)+2a\cdot \nabla_x u(x,y)\right )=0$ and therefore $\mathcal L^*$ is the operator 
$$y^\alpha\left(\Delta_x +2a\cdot \nabla_xD_y+2c\frac {a\cdot \nabla_x u}{y}+ B_y\right)$$ with the above oblique condition at $y=0$. 

\begin{prop}\label{generation L2}
 $\mathcal L$ and $\mathcal L^*$  generate contractive analytic semigroups $e^{z \mathcal L}$, $e^{z \mathcal L^*}$, $z\in\Sigma_{\frac{\pi}{2}-\arctan \frac{|a|}{1-|a|}}$,  in $L^2_{c-\alpha}$. Moreover the semigroups $(e^{t\mathcal L})_{t \geq 0}, (e^{t\mathcal L^*})_{t \geq 0}$ are positive and $L^p_{c-\alpha}$-contractive for $1 \leq p \leq \infty$.
\end{prop}
{\sc Proof.} We argue only for $\mathcal L$.
The  generation  result immediately follows from   Proposition \ref{prop-form} and \cite[Theorem 1.52]{Ouhabaz}.
The positivity follows by \cite[Theorem 2.6]{Ouhabaz} after observing that, if $u\in H^1_{\alpha,c}$, $u$ real, then $u^+\in H^1_{\alpha,c}$ and
\begin{align*}
	\mathfrak{a}(u^+,u^-)
	&:=
	\int_{\R^{N+1}_+} \langle \nabla u^+, \nabla u^-\rangle\,y^{c} dx\,dy+2\int_{\R^{N+1}_+}D_yu^+ a\cdot \nabla_xu^-\,y^{c} dx\,dy=0.
\end{align*}
Finally, the $L^\infty$-contractivity follows by \cite[Corollary 2.17]{Ouhabaz}  after observing that if $0\leq u\in H^1_{\alpha,c}$,  then $1\wedge u, (u-1)^+\in H^1_{\alpha,c}$ and, since $\nabla (1\wedge u)=\chi_{\{u<1\}}\nabla u$ 
and 
$\nabla (u-1)^+=\chi_{\{u>1\}}\nabla u$, one has 
\begin{align*}
	 \mathfrak{a}(1\wedge u,(u-1)^+)&=0.
\end{align*}
\qed

The Stein interpolation theorem then shows that the above semigroups are analytic in $L^p_{c-\alpha}$ for $1<p<\infty$, see \cite[Proposition 3.12]{Ouhabaz} and a result by Lamberton  yields maximal regularity in the same range, see \cite[Theorem 5.6]{KunstWeis}. Since our results are more general, we do not state these consequences here.

Our aim is to characterize the  domain of $\mathcal{L}$ in $L^2_{c-\alpha}$. As in \cite[Section 7.1]{MNS-Singular-Half-Space} and  \cite[Section 6.1]{MNS-CompleteDegenerate}, we can prove the following result.

\begin{teo}\label{Neumnan comp 2}
If  $c+1>|\alpha|$  then 
\begin{align*}
D(\mathcal L)&=W^{2,2}_\mathcal{N}(\alpha,\alpha,c-\alpha).
\end{align*}
In particular the set $C_c^\infty (\R^{N})\otimes \mathcal D$, see \eqref{defDcore}, 
is a core for $\mathcal L$ in $L^2_{c-\alpha}(\R^{N+1}_+)$.
\end{teo}
{\sc{Proof.}} The proof follows as in \cite[Proposition 7.3, Theorem 7.4]{MNS-Singular-Half-Space}  using  the  boundedness of the multipliers $|\xi|^2y^\alpha R_1(\xi)$, $\xi y^\alpha D_yR_1(\xi)$ in $L^2(\R^N; L^2_{c-\alpha}(\R_+))=L^2_{c-\alpha}$, proved in Proposition \ref{Mikhlin 0order}. Note that the condition $c+1>|\alpha|$ is that of the quoted proposition with $p=2$ and $m=c-\alpha$.\qed

\subsection{The operator $\mathcal L$ in $L^p_m$}

In this section we prove  domain characterization and maximal regularity for  $\mathcal L$ 
in $L^p_m$. For clarity reasons we often write $ \mathcal{L}_{m,p}$ to emphasize the underlying space on which the operator acts.

We shall use extensively the set  (finite sums below)
	$$C_c^\infty (\R^{N})\otimes\mathcal D=\left\{u(x,y)=\sum_i u_i(x)v_i(y), \  u_i \in C_c^\infty (\R^N), \  v_i \in \cal D \right \},$$
 where $\mathcal{D}$ is defined in \eqref{defD}. We refer to Appendix \ref{Weighted} for  further details and  note that 
  $\mathcal L$ is well defined on $C_c^\infty (\R^{N})\otimes\mathcal D \subset L^p_m$ when $(m+1)/p >\alpha^-$. 

The next results can be proved as in \cite[Lemma 7.5, Lemma 7.6, Theorem 7.7, Corollary 7.8]{MNS-Singular-Half-Space}.

\begin{lem}\label{lemma multipl}
Let $\alpha^-<\frac{m+1}p<c+1-\alpha$.
Then for any $\lambda\in \C^+$ the operators  
$$(\lambda- \mathcal{L}_{c-\alpha,2})^{-1},\quad y^\alpha \Delta_x (\lambda- \mathcal{L}_{c-\alpha,2})^{-1},\quad y^\alpha \nabla_xD_y (\lambda- \mathcal{L}_{c-\alpha,2})^{-1},\quad y^\alpha B^n_y (\lambda- \mathcal{L}_{c-\alpha,2})^{-1}$$ initially defined on $L^p_m\cap L^2_{c-\alpha}$ by Theorem \ref{Neumnan comp 2}, extend to bounded operators on $L^p_m$ which we denote respectively  by $\mathcal{R}(\lambda)$, $y^\alpha \Delta_x \mathcal{R}(\lambda)$, $y^\alpha \nabla_xD_y \mathcal{R}(\lambda)$, $y^\alpha B^n_y \mathcal{R}(\lambda)$. Moreover  the family  $\left\{\lambda \mathcal{R}(\lambda):\lambda\in\C^+\right\}$ is $\mathcal{R}$-bounded on $L^p_m$.
\end{lem}

\begin{prop} \label{generazione} If  $ \alpha^-<\frac{m+1}p<c+1-\alpha$, an extension $\mathcal L_{m,p}$ of the operator $\mathcal L$, initially defined on $C_c^\infty (\R^{N})\otimes  \mathcal{D}$,  generates a bounded analytic semigroup  in $L^p_m(\R_+^{N+1})$ which has maximal regularity and it is consistent with the semigroup generated by $\mathcal L_{c-\alpha,2}$ in $L^2_{c-\alpha}(\R^{N+1}_+)$. 
\end{prop}

Finally we characterize the domain of $\mathcal L_{m,p}$.
\begin{teo} \label{generazione1} 
If $\alpha^-<\frac{m+1}p<c+1-\alpha$, then
\begin{align*}
D(\mathcal L_{m,p})&=W^{2,p}_{\mathcal{N}}(\alpha,\alpha,m)
\end{align*}
 and in particular $C_c^\infty (\R^{N})\otimes \mathcal D$ is a core for $\mathcal L_{m,p}$.
\end{teo}

\begin{cor} \label{omogeneo}
Under the hypotheses of Theorem \ref{generazione1} we have for every $u \in W^{2,p}_{\mathcal{N}}(\alpha,\alpha,m)$
$$
\|y^\alpha D_{x_i x_j} u\|_{L^p_{m}} +\|y^\alpha D_{yy} u\|_{L^p_{m}}+\|y^\alpha D_{x_iy} u\|_{L^p_{m}}+\|y^{\alpha-1} D_{y} u\|_{L^p_{m}}\leq C\| \mathcal Lu\|_{L^p_{m}}.
$$
\end{cor}

\section{Consequences for more general operators}\label{Consequenses}

The isometry  introduced in Section \ref{Weighted} allows to deduce generation and domain properties in  $L^p_m$ for more general operators  of the form
\begin{equation*}
	\mathcal L =y^{\alpha_1}\Delta_x+2y^\frac{\alpha_1+\alpha_2}{2}a\cdot\nabla_xD_y+y^{\alpha_2}\left(D_{yy}+\frac{c}{y}D_y\right),
\end{equation*}
with $\alpha_1\alpha_2\in\R$,  $\alpha_2<2$, $\alpha_2-\alpha_1<2$.

\begin{teo} \label{complete-Bessel} 
Let $\alpha_2<2$,  $\alpha_2-\alpha_1<2$ and
$$\alpha_1^{-} <\frac{m+1}p<c+1-\alpha_2.$$   
  Then $\mathcal L$ with domain $ W^{2,p}_{\mathcal{N}}\left(\alpha_1,\alpha_2,m\right)$ generates a bounded analytic semigroup  in $L^p_m$ which has maximal regularity.
\end{teo}
{\sc Proof.} We use the isometry 
$$T_{\frac{\alpha_1-\alpha_2}{2}}:L^p_{\tilde{m}}\to L^p_m,\qquad \tilde m=\frac{2m-\alpha_1+\alpha_2}{\alpha_1-\alpha_2+2}$$
which, according to Proposition \ref{Isometry action},   transforms  $\mathcal L$ into  
$$
T_{\frac{\alpha_1-\alpha_2}{2}}^{-1} \,\mathcal L\, T_{\frac{\alpha_1-\alpha_2}{2}}=y^\alpha\Delta_x+y^\alpha 2\left(\frac{\alpha_1-\alpha_2+2}{2}\right)a\cdot\nabla_xD_y+\left(\frac{\alpha_1-\alpha_2+2}{2}\right)^2y^\alpha \tilde B_y
$$ where $$\alpha=\frac{2\alpha_1}{\alpha_1-\alpha_2+2}, \quad \tilde B_y=D_{yy}+\frac{\tilde c}{y}D_y ,\quad
 \tilde c=\frac{2c+\alpha_1-\alpha_2}{\alpha_1-\alpha_2+2}.$$
Observe the assumptions on the parameters translates into  $\alpha<2$ and $\alpha^- < \frac{\tilde m +1}{p} < \tilde c+1-\alpha$. The generation properties and maximal regularity  for $\mathcal L$ in $L^p_m$ are then  immediate consequence of the same properties for the operator studied before.
Concerning the domain, we have
$$D(\mathcal L)=T_{\frac{\alpha_1-\alpha_2}{2}}\left( W^{2,p}_{\mathcal{N}}\left (\alpha, \alpha,  \tilde m \right)\right)$$
which, by Proposition \ref{Isometry action der}, coincides with $W^{2,p}_{\mathcal{N}}\left(\alpha_1,\alpha_2,m\right)$.
\qed


\medskip
Results for more general operators  follow by linear change of variables, as we explain below.  We  consider the operator   in $\R^{N+1}_+$
\begin{align*}
	\mathcal L&=y^{\alpha_1}\mbox{Tr }\left(QD^2_xu\right)+2y^{\frac{\alpha_1+\alpha_2}{2}}q\cdot \nabla_xD_y+y^{\alpha_2}\left(\gamma D_{yy}+\frac{c}{y}D_y\right)\\[1ex]
	 &=y^{\alpha_1}\sum_{i,j=1}^{N}q_{ij}D_{x_ix_j}+2y^\frac{\alpha_1+\alpha_2}{2}\sum_{i=1}^Nq_iD_{x_iy}+y^{\alpha_2}\left(\gamma D_{yy}+\frac{c}{y}D_y\right).
\end{align*}
Here $Q$ is the $N\times N$ matrix $(q_{ij})$,  $q=(q_1, \dots, q_N)$ and  we assume that the quadratic form $Q(\xi,\xi)+\gamma \eta^2+2 q\cdot\xi \,\eta$ is positive definite.
Through a linear change of variables in the $x$ variables the term $\sum_{i,j=1}^{N}q_{ij}D_{x_ix_j}$ is transformed into $\gamma\Delta_x$ and all the results of Section \ref{Sect DOm max} hold, replacing $c$ with $ \frac c\gamma$ in the statements (the condition $|a|<1$ of Section \ref{Sect DOm max} is satisfied since the change of variables preserves the ellipticity). 
The case of variable coefficients can also be handled by freezing the coefficients and will be done in the future to deal with degenerate problems in bounded domains.

\bigskip

We can therefore deduce  results also for the last operator whose proofs follow directly from   Theorem \ref{complete-Bessel}.

\begin{teo} \label{complete-oblique}
Let  {$c\in\R$ and} $
\left(
\begin{array}{c|c}
	 Q  & { q}^t \\[1ex] \hline
	 q& \gamma
\end{array}\right)$
an  elliptic matrix and  	let $\alpha_1,\alpha_2\in\R$ such that $\alpha_2<2$,  $\alpha_2-\alpha_1<2$ and
	$$\alpha_1^{-} <\frac{m+1}p<\frac{c}{\gamma}+1-\alpha_2.$$   
	Then the operator
	\begin{align*}
	\mathcal L=y^{\alpha_1}\mbox{Tr }\left(QD^2_xu\right)+2y^{\frac{\alpha_1+\alpha_2}{2}}q\cdot \nabla_xD_y+\gamma y^{\alpha_2} D_{yy}+{c y^{\alpha_2-1}D_y}
	\end{align*} endowed with domain {$W^{2,p}_{\mathcal N}\left(\alpha_1,\alpha_2,m\right)$}  generates a bounded analytic semigroup  in $L^p_m$ which has maximal regularity.
\end{teo}

{We refer the reader to \cite{Negro-AlphaDirichlet}  for a further generalization of $\mathcal L$ involving Dirichlet or oblique derivative boundary conditions.}

\section{Appendix A: Vector-valued harmonic analysis}
%

We review some results on vector-valued multiplier theorems referring the reader to \cite{DenkHieberPruss}, \cite{Pruss-Simonett} or \cite{KunstWeis} for all proofs.

Let ${\cal S}$ be a subset of $B(X)$, the space of all bounded linear operators on a Banach space $X$. ${\cal S}$ is $\mathcal R$-bounded if there is a constant $C$ such that
$$
\|\sum_i \eps_i S_i x_i\|_{L^p(\Omega; X)} \le C\|\sum_i \eps_i  x_i\|_{L^p(\Omega; X)} 
$$
for every finite sum as above, where $(x_i ) \subset X, (S_i) \subset {\cal S}$ and $\eps_i:\Omega \to \{-1,1\}$ are independent and symmetric random variables on a probability space $\Omega$. In particular ${\cal S}$ is a bounded subset of $B(X)$. The smallest constant $C$ for which the above definition holds is the $\mathcal R$-bound of $\mathcal S$,  denoted by $\mathcal R(\mathcal S)$.
It is well-known that this definition does not  depend on $1 \le p<\infty$ (however, the constant $\mathcal R(\mathcal S)$ does) and that $\mathcal R$-boundedness  is equivalent to boundedness when $X$ is an Hilbert space.
When $X$ is an $L^p(\Sigma)$ space (with respect to any $\sigma$-finite measure defined on a $\sigma$-algebra $\Sigma$), testing  $\mathcal R$-boundedness is equivalent to proving square functions estimates, see \cite[Remark 2.9]{KunstWeis}.

\begin{prop}\label{Square funct R-bound} Let ${\cal S} \subset B(L^p(\Sigma))$, $1<p<\infty$. Then ${\cal S}$ is $\mathcal R$-bounded if and only if there is a constant $C>0$ such that for every finite family $(f_i)\in L^p(\Sigma), (S_i) \in {\cal S}$
	$$
	\left\|\left (\sum_i |S_if_i|^2\right )^{\frac{1}{2}}\right\|_{L^p(\Sigma)} \le C\left\|\left (\sum_i |f_i|^2\right)^{\frac{1}{2}}\right\|_{L^p(\Sigma)}.
	$$
\end{prop}
The best constant $C$ for which the above square functions estimates hold satisfies $\kappa^{-1} C \leq \mathcal R(\mathcal S) \leq \kappa C$ for a suitable $\kappa>0$ (depending only on $p$).  Using the proposition above, $\mathcal R$-boundedness follows from domination by a positive $\mathcal R$-bounded family.
\begin{cor} \label{domination}
	Let  ${\cal S}, {\cal T} \subset B(L^p(\Sigma))$, $1<p<\infty$ and assume that $\cal T$ is an $\mathcal R$ bounded family of positive operators and that for every $S \in \cal S$ there exists $T \in \cal T$ such that $|Sf| \leq T|f|$ pointwise, for every $f \in L^p(\Sigma)$. Then ${\cal S}$ is $\mathcal R$-bounded.
\end{cor}

We also need the following result about the integral mean of a $\mathcal R$-bounded  family of operator which we state in the version we use.
\begin{prop}{\em \cite[Corollary 2.14]{KunstWeis}} \label{Mean R-bound}
	Let  $X$ be a Banach space and let  ${\cal F}\subset B(X)$ be an $\mathcal R$-bounded  family of operator. For every strongly measurable $N:\Sigma\to B(X)$ on a $\sigma$-finite measure space  $(\Sigma,\mu)$ with values in 
	$\cal F$ and every $h\in L^1\left(\Sigma,\mu \right)$ we define the operator  $T_{N,\cal F}\in B(X)$ by
	\begin{align*}
		T_{N,\cal F} x=\int_{\Sigma} h(\omega)N(\omega)xd\mu(\omega),\qquad x\in X.
	\end{align*}
	Then the family 
	\begin{align*}
		\mathcal C=\left\{T_{N,\cal F}\;:\; \|h\|_{L^1}\leq 1, N\;\text{as above}\right\}
	\end{align*} is $\mathcal R$ bounded and $\mathcal R(\mathcal C)\leq 2\mathcal R(\mathcal F)$.
\end{prop}
Let $(A, D(A))$ be a sectorial operator in a Banach space $X$; this means that $\rho (-A) \supset \Sigma_{\pi-\phi}$ for some $\phi <\pi$ and that $\lambda (\lambda+A)^{-1}$ is bounded in $\Sigma_{\pi-\phi}$. The infimum of all such $\phi$ is called the spectral angle of $A$ and denoted by $\phi_A$. Note that $-A$ generates an analytic semigroup if and only if $\phi_A<\pi/2$. The definition of $\mathcal  R$-sectorial operator is similar, substituting boundedness of $\lambda(\lambda+A)^{-1}$ with $\mathcal R$-boundedness in $\Sigma_{\pi-\phi}$. As above one denotes by $\phi^R_A$ the infimum of all $\phi$ for which this happens; since $\mathcal R$-boundedness implies boundedness, we have $\phi_A \le \phi^R_A$.

\medskip

The $\mathcal R$-boundedness of the resolvent characterizes the regularity of the associated inhomogeneous parabolic problem, as we explain now.

An analytic semigroup $(e^{-tA})_{t \ge0}$ on a Banach space $X$ with generator $-A$ has
{\it maximal regularity of type $L^q$} ($1<q<\infty$)
if for each $f\in L^q([0,T];X)$ the function
$t\mapsto u(t)=\int_0^te^{-(t-s)A})f(s)\,ds$ belongs to
$W^{1,q}([0,T];X)\cap L^q([0,T];D(A))$.
This means that the mild solution of the evolution equation
$$u'(t)+Au(t)=f(t), \quad t>0, \qquad u(0)=0,$$
is in fact a strong solution and has the best regularity one can expect.
It is known that this property does not depend on $1<q<\infty$ and $T>0$.
A characterization of maximal regularity is available in UMD Banach spaces, through the $\mathcal  R$-boundedness of the resolvent in a suitable sector $\omega+\Sigma_{\phi}$, with $\omega \in \R$ and $\phi>\pi/2$ or, equivalently, of the scaled semigroup $e^{-(A+\omega')t}$ in a sector around the positive axis. In the case of $L^p$ spaces it can be restated in the following form,  see \cite[Theorem 1.11]{KunstWeis}

\begin{teo}\label{MR} Let $(e^{-tA})_{t \ge0}$ be a bounded analytic semigroup in $L^p(\Sigma)$, $1<p<\infty$,  with generator $-A$. Then $T(\cdot)$ has maximal regularity of type $L^q$  if and only if the set $\{\lambda(\lambda+A)^{-1}, \lambda \in  \Sigma_{\pi/2+\phi} \}$ is $\mathcal R$- bounded for some $\phi>0$. In an equivalent way, if and only if 
	there are constants $0<\phi<\pi/2 $, $C>0$ such that for every finite sequence $(\lambda_i) \subset \Sigma_{\pi/2+\phi}$, $(f_i) \subset  L^p$
	$$
	\left\|\left (\sum_i |\lambda_i (\lambda_i+A)^{-1}f_i|^2\right )^{\frac{1}{2}}\right\|_{L^p(\Sigma)} \le C\left\|\left (\sum_i |f_i|^2\right)^{\frac{1}{2}}\right\|_{L^p(\Sigma)}
	$$
	or, equivalently, there are constants $0<\phi'<\pi/2 $, $C'>0$ such that  for every finite sequence
	$(z_i) \subset \Sigma_{\phi'}$, $(f_i) \subset  L^p$
	$$
	\left\|\left (\sum_i |e^{-z_i A}f_i|^2\right )^{\frac{1}{2}}\right\|_{L^p(\Sigma)} \le C'\left\|\left (\sum_i |f_i|^2\right)^{\frac{1}{2}}\right\|_{L^p(\Sigma)}.
	$$
\end{teo}
\medskip

Finally we state  a version of the operator-valued Mikhlin multiplier theorem in the N-dimensional case, see e.g. \cite[Corollary 8.3.22]{WeisBook2}.

\begin{teo}   \label{marcinkiewicz}
	Let $1<p<\infty$, $M\in C^N(\R^N\setminus \{0\}; B(L^p(\Sigma))$ be such that  the set
	$$\left \{\xi^{\alpha}D^\alpha_\xi M(\xi): \xi\in \R^{N}\setminus\{0\}, \ \alpha\in\{0,1\}^N \right \}$$
	is $\mathcal{R}$-bounded.
	Then the operator $T_M={\cal F}^{-1}M {\cal F}$ is bounded in $L^p(\R^N, L^p(\Sigma))$, where $\cal{F}$ denotes the Fourier transform.
\end{teo}

\section{Appendix B: Weighted spaces and similarity transformations} \label{Weighted} 
Let $p>1$, $m, \alpha_1,\alpha_2 \in \R$ such that
\begin{align*}
	\alpha_2<2,\qquad \alpha_2-\alpha_1<2,\qquad \alpha_1^{-} <\frac{m+1}p.
\end{align*}
In order to describe the domain of  the operator \begin{align*}
	y^{\alpha_1}\Delta_x+2y^\frac{\alpha_1+\alpha_2}{2}a\cdot\nabla_xD_y+y^{\alpha_2}\left(D_{yy}+\frac{c}{y}D_y\right),
\end{align*}
we collect in this section  the main results concerning anisotropic weighted Sobolev spaces, referring to \cite{MNS-Sobolev} for further details  and all the relative proofs.
We define the Sobolev space
\begin{align*}
	W^{2,p}(\alpha_1,\alpha_2,m)=&\left\{u\in W^{2,p}_{loc}(\R^{N+1}_+):\ u,\  y^{\alpha_1} D_{x_ix_j}u,\ y^\frac{\alpha_1}{2} D_{x_i}u,  \right.\\[1ex]
	&\left.\hspace{20ex}y^{\alpha_2}D_{yy}u,\ y^{\frac{\alpha_2}{2}}D_{y}u,\,y^\frac{\alpha_1+\alpha_2}{2} D_{y}\nabla_x u\in L^p_m\right\}
\end{align*}
which is a Banach space equipped with the norm
\begin{align*}
	\|u\|_{W^{2,p}(\alpha_1,\alpha_2,m)}=&\|u\|_{L^p_m}+\sum_{i,j=1}^n\|y^{\alpha_1} D_{x_ix_j}u\|_{L^p_m}+\sum_{i=1}^n\|y^{\frac{\alpha_1}2} D_{x_i}u\|_{L^p_m}\\[1ex]
	&+\|y^{\alpha_2}D_{yy}u\|_{L^p_m}+\|y^{\frac{\alpha_2}{2}}D_{y}u\|_{L^p_m}+\|y^\frac{\alpha_1+\alpha_2}{2} D_{y}\nabla_x u\|_{L^p_m}.
\end{align*}
Next we add a Neumann boundary condition for $y=0$  in the form $y^{\alpha_2-1}D_yu\in L^p_m$ and set
\begin{align*}
	W^{2,p}_{\mathcal N}(\alpha_1,\alpha_2,m)=\{u \in W^{2,p}(\alpha_1,\alpha_2,m):\  y^{\alpha_2-1}D_yu\ \in L^p_m\}
\end{align*}
with the norm
$$
\|u\|_{W^{2,p}_{\mathcal N}(\alpha_1,\alpha_2,m)}=\|u\|_{W^{2,p}(\alpha_1,\alpha_2,m)}+\|y^{\alpha_2-1}D_yu\|_{ L^p_m}.
$$

\begin{os}\label{Os Sob 1-d}
	With obvious changes we consider also the analogous Sobolev spaces $W^{2,p}(\alpha,m)$ and $W^{2,p}_{\cal N}(\alpha, m)$ on $\R_+$. 
	For example  we have
	$$W^{2,p}_{\mathcal N}(\alpha,m)=\left\{u\in W^{2,p}_{loc}(\R_+):\ u,\    y^{\alpha}D_{yy}u,\ y^{\frac{\alpha}{2}}D_{y}u,\ y^{\alpha-1}D_{y}u\in L^p_m\right\}.$$
	All the results of this section will be valid also in $\R_+$ changing (when it appears) the condition $\alpha_1^{-} <\frac{m+1}p$  to $0<\frac{m+1}p$.
\end{os}

The next result clarifies in which sense the condition $y^{\alpha_2-1}D_y u \in L^p_m$ is a Neumann boundary condition.

\begin{prop}{\em \cite[Proposition 4.3]{MNS-Sobolev}} \label{neumann} The following assertions hold.
	\begin{itemize} 
		\item[(i)] If $\frac{m+1}{p} >1-\alpha_2$, then $W^{2,p}_{\mathcal N}(\alpha_1, \alpha_2, m)=W^{2,p}(\alpha_1, \alpha_2, m)$.
		\item[(ii)] If $\frac{m+1}{p} <1-\alpha_2$, then $$W^{2,p}_{\mathcal N}(\alpha_1, \alpha_2, m)=\{u \in W^{2,p}(\alpha_1, \alpha_2, m): \lim_{y \to 0}D_yu(x,y)=0\ {\rm for\ a.e.\   x \in \R^N }\}.$$
	\end{itemize}
	In both cases (i) and (ii), the norm of $W^{2,p}_{\mathcal N}(\alpha_1, \alpha_2, m)$ is equivalent to that of $W^{2,p}(\alpha_1, \alpha_2, m)$.
\end{prop}

The next results  show the density  of smooth functions in $W^{2,p}_{\mathcal N}(\alpha_1,\alpha_2,m)$. Let
\begin{equation} \label{defC}
	\mathcal{C}:=\left \{u \in C_c^\infty \left(\R^N\times[0, \infty)\right), \ D_y u(x,y)=0\  {\rm for} \ y \leq \delta\ {\rm  and \ some\ } \delta>0\right \},
\end{equation}
its one dimensional version 
\begin{equation} \label {defDcore}
	\mathcal{D}=\left \{u \in C_c^\infty ([0, \infty)), \ D_y u(y)=0\  {\rm for} \ y \leq \delta\ {\rm  and \ some\ } \delta>0\right \}
\end{equation}
and finally (finite sums below)
$$C_c^\infty (\R^{N})\otimes\mathcal D=\left\{u(x,y)=\sum_i u_i(x)v_i(y), \  u_i \in C_c^\infty (\R^N), \  v_i \in \cal D \right \}\subset \mathcal C.$$
\begin{teo} {\em \cite[Theorem 4.9]{MNS-Sobolev}}\label{core gen}
	$C_c^\infty (\R^{N})\otimes\mathcal D$ is dense in $W^{2,p}_{\mathcal N}(\alpha_1,\alpha_2,m)$.
\end{teo}

Note that the condition $(m+1)/p>\alpha_1^-$, or $m+1>0$ and $(m+1)/p+\alpha_1>0$, is necessary for the inclusion  $C_c^\infty (\R^{N})\otimes\mathcal D \subset W^{2,p}_{\mathcal N}(\alpha_1,\alpha_2,m)$. 
\medskip

In the 1-dimensional case we need also the following density result proved in  \cite[Theorem 8.5]{MNS-PerturbedBessel}.
\begin{teo} \label{core1-d}  Let  $\alpha<2$, $\mu>0$, $c\in\R$.  Then for any $1<p<\infty$ such that $\alpha^-<\frac{m+1}p<c+1-\alpha$, the set $\mathcal {D}$ 
is dense in  $W^{2,p}_{\mathcal N}(\alpha,m)\cap L^p_{m+\alpha p}$.
\end{teo}

We consider now, for $\beta \in\R$, $\beta\neq -1$, the transformation
\begin{align}\label{Gen Kelvin def}
	T_{\beta\,}u(x,y)&:=|\beta+1|^{\frac 1 p}u(x,y^{\beta+1}),\quad (x,y)\in\R^{N+1}_+.
\end{align}
Observe that
$$ T_{\beta\,}^{-1}=T_{-\frac{\beta}{\beta+1}\,}.$$

\begin{prop}\label{Isometry action der} Let $1\leq p\leq \infty$, $\beta \in\R$, $\beta\neq -1$ and $m\in\R$. The following properties hold.
	\begin{itemize}
		\item[(i)] $T_{\beta\,}$ maps isometrically  $L^p_{\tilde m}$ onto $L^p_m$  where $ \tilde m=\frac{m-\beta}{\beta+1}$.
		\item[(ii)] $W^{2,p}_{\mathcal{N}}(\alpha_1,\alpha_2,m)=T_{\beta}\left(W^{2,p}_{\mathcal{N}}(\tilde \alpha_1,\tilde\alpha_2,\tilde m)\right),\qquad \tilde\alpha_1=\frac{\alpha_1}{\beta+1},\quad \tilde\alpha_2=\frac{\alpha_2+2\beta}{\beta+1}.$
	\end{itemize}
	In particular choosing $\beta=\frac{\alpha_1-\alpha_2}2$ and setting $\tilde \alpha=\frac{2\alpha_1}{\alpha_1-\alpha_2+2}$ one has
	\begin{align*}
		W^{2,p}_{\mathcal{N}}(\alpha_1,\alpha_2,m)=T_{\frac{\alpha_1-\alpha_2}2}\left(W^{2,p}_{\mathcal{N}}(\tilde \alpha,\tilde \alpha ,\tilde m)\right),\qquad \tilde \alpha=\frac{2\alpha_1}{\alpha_1-\alpha_2+2},\quad \tilde m=\frac{2m-\alpha_1+\alpha_2}{\alpha_1-\alpha_2+2}.
	\end{align*}
\end{prop}
{\sc{Proof.}} See \cite[Lemma 2.1, Proposition 2.2]{MNS-Sobolev} with $k=0$ and 	$T_{\beta}=T_{0,\beta}$.\qed

\begin{os}
	It is essential to deal with  $W^{2,p}_{\mathcal{N}}(\alpha_1,\alpha_2,m)$: in  general the map $T_{\beta}$ does  not  transform  $W^{2,p}(\tilde \alpha_1,\tilde\alpha_2,\tilde m)$ into $W^{2,p}(\alpha_1,\alpha_2,m)$.
\end{os}

We consider now  the  operators
$$\mathcal L=y^{\alpha_1}\Delta_x+2y^{\frac{\alpha_1+\alpha_2}2}\left(a,\nabla_xD_y\right)+y^{\alpha_2} B_y,\qquad a \in\R^N,\ |a|<1$$ 
in the space $L^p_m=L^p_m(\R^{N+1}_+)$. Here $B$ is the Bessel operator
$$
B=D_{yy}+\frac{c}{y}D_y,\qquad c>-1
$$ on the  half line $\R_+=]0, \infty[$ (often we write $B_y$ to indicate that it acts with respect to the $y$ variable). The condition $|a| <1$ is equivalent to the ellipticity of the top order coefficients
We investigate when these operators can be transformed one into the other by means of the transformation \eqref{Gen Kelvin def}.

\begin{prop}\label{Isometry action}  Let 
	$T_{\beta\,}$ be the isometry defined in \eqref{Gen Kelvin def}.   Then for every  $u\in W^{2,1}_{loc}\left(\R^{N+1}_+\right)$ one has
	\begin{align*}
		&T_{\beta\,}^{-1} \Big(y^{\alpha_1}\Delta_x+2y^{\frac{\alpha_1+\alpha_2}2}\left(a,\nabla_xD_y\right)+y^{\alpha_2} B_y\Big)T_{\beta\,}u\\[1ex]
		&=\Big(y^{\frac{\alpha_1}{\beta+1}}\Delta_x+2(\beta+1)y^{\frac{\alpha_1+\alpha_2+2\beta}{2(\beta+1)}}\left(a,\nabla_xD_y\right)+(\beta+1)^2y^{\frac{\alpha_2+2\beta}{\beta+1}}\tilde{ B}_y\Big) u
	\end{align*}
	where 
	\begin{align}\label{tilde b}
		\tilde B =D_{yy}+\frac{\tilde c}{y}D_y,\qquad 
		\tilde c&=\frac{c+\beta}{\beta+1}.
	\end{align}
	In particular choosing $\beta=\frac{\alpha_1-\alpha_2}2$ and setting $\tilde \alpha=\frac{2\alpha_1}{\alpha_1-\alpha_2+2}$ one has 
	\begin{align*}
		&T_{\beta\,}^{-1} \Big(y^{\alpha_1}\Delta_x+2y^{\frac{\alpha_1+\alpha_2}2}\left(a,\nabla_xD_y\right)+y^{\alpha_2} B_y\Big)T_{\beta\,}u\\[1ex]
		&\hspace{15ex}=y^{\tilde \alpha}\Big(\Delta_x+2(\beta+1)\left(a,\nabla_xD_y\right)+(\beta+1)^2\tilde{ B}_y\Big) u
	\end{align*}
\end{prop}
{\sc{Proof.}} The proof follows using \cite[Proposition 3.1, Proposition 3.2]{MNS-CompleteDegenerate} with $k=0$ and the equalities
\begin{align*}
	y^\alpha T_{\beta\,}u=T_{\beta\,}(y^{\frac{\alpha}{\beta+1}}u),\quad  D_{x_i}(T_{\beta\,}u)=T_{\beta}\left(D_{x_i} u\right),\quad 
	D_{xy} T_{\beta\,}u=T_{\beta\,}\left((\beta+1)y^{\frac{\beta}{\beta+1}}D_{xy}u\right).
\end{align*}
\qed

\bibliography{../../TexBibliografiaUnica/References}

\end{document}